\documentclass{amsart}
\usepackage[utf8]{inputenc}

\usepackage{amsmath}
\usepackage{amssymb,amsfonts}

\usepackage[mathscr]{euscript}
\usepackage[utf8]{inputenc}
\usepackage[pdftex]{graphicx}
\usepackage{xypic}
\usepackage{enumitem}
\usepackage{chngcntr}
\usepackage{txfonts}
\usepackage{calligra}
\usepackage[T1]{fontenc}
\usepackage{frcursive}
\usepackage[T1]{fontenc}

\usepackage{todonotes}
\usepackage[capitalize]{cleveref}
\usepackage{verbatim}
\usepackage[normalem]{ulem}

\theoremstyle{definition}

\theoremstyle{definition}

\newcommand{\To}{\Rightarrow}

\newcommand{\Aa}{\mathcal{A}}
\newcommand{\Bb}{\mathcal{B}}
\newcommand{\Cc}{\mathcal{C}}

\newcommand{\Ee}{\mathcal{E}}
\newcommand{\Ff}{\mathcal{F}}

\newcommand{\Hh}{\mathcal{H}}

\newcommand{\Jj}{\mathcal{J}}

\newcommand{\Mm}{\mathcal{M}}

\newcommand{\Rr}{\mathcal{R}}
\newcommand{\Ss}{\mathcal{S}}

\newcommand{\Vv}{\mathcal{V}}

\newcommand{\AAA}{\mathbb{A}}
\newcommand{\BB}{\mathbb{B}}
\newcommand{\CC}{\mathbb{C}}

\newcommand{\EE}{\mathbb{E}}
\newcommand{\FF}{\mathbb{F}}

\newcommand{\JJ}{\mathbb{J}}

\newcommand{\MM}{\mathbb{M}}
\newcommand{\NN}{\mathbb{N}}

\newcommand{\RR}{\mathbb{R}}
\newcommand{\SSS}{\mathbb{S}}

\newcommand{\UU}{\mathbb{U}}
\newcommand{\VV}{\mathbb{V}}

\newcommand{\asf}{\mathsf{a}}

\newcommand{\csf}{\mathsf{c}}

\newcommand{\lsf}{\mathsf{l}}

\newcommand{\rsf}{\mathsf{r}}

\newcommand{\Abf}{\mathbf{A}}
\newcommand{\Bbf}{\mathbf{B}}
\newcommand{\Cbf}{\mathbf{C}}

\newcommand{\un}{\sf 1}
\newcommand{\zero}{\sf 0}

\title{The groupoid of finite sets is biinitial in the 2-category of rig categories}

\author{Josep Elgueta}
\address{Departament de Matem\`atiques Universitat Polit\`ecnica de Catalunya}
\email{josep.elgueta@upc.edu}
\thanks{This work had financial support from the Generalitat de Catalunya (Project: 2014 SGR 634), and the Ministerio de Econom\'ia y Competitividad of Spain (Project: MTM2015-69135-9).}



\begin{document}

\begin{abstract}
The groupoid of finite sets has a ``canonical'' structure of a symmetric 2-rig with the sum and product respectively given by the coproduct and product of sets. This 2-rig $\widehat{\FF\SSS et}$ is just one of the many non-equivalent categorifications of the commutative rig $\NN$ of natural numbers, together with the rig $\NN$ itself viewed as a discrete rig category, the whole category of finite sets, the category of finite dimensional vector spaces over a field $k$, etc. In this paper it is shown that  $\widehat{\FF\SSS et}$ is the right categorification of $\NN$ in the sense that it is biinitial in the 2-category of rig categories, in the same way as $\NN$ is initial in the category of rigs. As a by-product, an explicit description of the homomorphisms of rig categories from a suitable version of $\widehat{\FF\SSS et}$ into any (semistrict) rig category $\SSS$ is obtained in terms of a sequence of automorphisms of the objects $1+\stackrel{n)}{\cdots}+1$ in $\SSS$ for each $n\geq 0$.
\end{abstract}

\maketitle

\section{Introduction}

A {\em rig} (a.k.a. a semiring) is a ring without {\em n}egatives, i.e. an (additive) abelian monoid $(S,+,0)$ equipped with an additional (multiplicative) monoid structure $(\,\bullet\,,1)$ such that $\bullet$ distributes over $+$ on both sides, and $0\bullet x=x\bullet 0=0$ for each $x\in S$. A paradigmatic example is the set $\NN$ of nonnegative integers with the usual sum and product. The rig is called {\em commutative} when $\bullet$ is abelian; for instance, $(\NN,+,\bullet\,,0,1)$ is a commutative rig. Rigs (resp. commutative rigs) are the objects of a category $\Rr ig$ (resp. $\Cc\Rr ig$) having as morphisms the maps that preserve both $+$ and $\bullet$\,, and the corresponding neutral elements.

We are interested in the categorical analog of a (commutative) rig. It is usually called a ({\em symmetric}) {\em rig category}, or a ({\em symmetric}) {\em bimonoidal category}. The last name, however, is confusing because a bimonoid in the set-theoretic context is a set simultaneously equipped with compatible monoid and comonoid structures, and not a set with two monoid structures one of them distributing over the other. 

The precise definition of a (symmetric) rig category is due to Laplaza \cite{Laplaza-1972}, and goes back to the 1970s (see also \cite{Kelly-1974}). Roughly, it is a category $\Ss$ equipped with functorial operations analogous to those of a rig, and satisfying all rig axioms up to suitable natural isomorphisms. More precisely, $\Ss$ must be equipped with an (additive) symmetric monoidal structure $(+,0,a,c,l,r)$, a (multiplicative) monoidal structure $(\bullet,1,a',l',r')$ (including a muliplicative commutator $c'$ in case the rig category is symmetric), and distributor and absorbing natural isomorphisms
\begin{align*}
d_{x,y,z}&: x\bullet(y+z)\stackrel{\cong}{\to} x\bullet y+x\bullet z
\\
d'_{x,y,z}&:(x+y)\bullet z\stackrel{\cong}{\to} x\bullet z+y\bullet z
\\
n_x&:x\bullet 0\stackrel{\cong}{\to} 0
\\
n'_x&:0\bullet x\stackrel{\cong}{\to} 0
\end{align*}
making commutative the appropriate `coherence diagrams'. For short, we shall denote by $\SSS$ the whole data defining a (symmetric) rig category. A paradigmatic example is the category $\Ff\Ss et$ of finite sets and maps between them, with $+$ and $\bullet$ respectively given by the disjoint union and the cartesian product of finite sets. The rig category so defined $\FF\SSS et$ is symmetric with $c'$ given by the canonical isomorphisms of sets $S\times T\cong T\times S$. Of course, rig categories (resp. symmetric rig categories) are the objects of a {\em 2-category} $\mathbf{RigCat}$ (resp. $\mathbf{SRigCat}$) whose 1- and 2-cells are given by the appropiate type of functors and natural transformations between these. The precise definitions are given in Section~\ref{categories_rig}.

$\FF\SSS et$ is not just a symmetric rig category. It is a {\em categorification} of the commutative rig $\NN$, in the sense that the set of isomorphisms classes of objects in $\Ff\Ss et$ with the rig structure induced by $+$ and $\bullet$ is isomorphic to $\NN$ with its canonical rig structure; for more on the idea of categorification, see \cite{Baez-Dolan-1998}. In fact, there are many other non-equivalent categorifications of $\NN$ as a rig, such as the symmetric rig category $\FF\VV ect_k$ of finite dimensional vector spaces over any given field $k$, with $+$ and $\bullet$ respectively given by the direct sum and the tensor product of vector spaces, or the rig $\NN$ itself viewed as a discrete category with only the identity morphisms. What is then the right categorical analog of $\NN$ as a rig? Of course, the answer depends on what we mean by the ``right categorical analog'' of $\NN$ as a rig. As it is well known, $\NN$ is an {\em initial} object in the category of rigs, i.e. for every rig $(S,+,\bullet\,,0,1)$ there is one and only one rig homomorphism from $\NN$ to $S$, namely, the map $\varphi:\NN\to S$ given by $\varphi(0)=0$ and $\varphi(n)=1+\stackrel{n)}{\cdots}+1$ for each $n\geq 1$. Hence it is reasonable to look for the symmetric rig category having the analogous categorical property.

It is a conjecture, apparently due to John Baez,~\footnote{See the nLab webpage {\em https://ncatlab.org/nlab/show/rig+category}.} that the right categorical analog of $\NN$ in this sense is the {\em groupoid} of finite sets and the bijections between them equipped with the symmetric rig category structure inherited from $\FF\SSS et$. Thus if we denote by $\widehat{\FF\SSS et}$ this symmetric rig category, Baez's conjecture is that $\widehat{\FF\SSS et}$ is {\em biinitial} in $\mathbf{SRigCat}$, i.e. for every symmetric rig category $\SSS$ the category of symmetric rig category homomorphisms of symmetric rig categories from $\widehat{\FF\SSS et}$ to $\SSS$ is expected to be equivalent to the terminal category with only one object and its identity morphism. The purpose of this paper is to prove that $\widehat{\FF\SSS et}$ is in fact biinitial in $\mathbf{RigCat}$, in the same way as $\NN$ is initial in the category of rigs. For instance, up to isomorphism, the usual free vector space functor from the groupoid $\widehat{\Ff\Ss et}$ to $\Ff\Vv ect_k$ is the unique functor which extends to a homomorphism of rig categories $\widehat{\FF\SSS et}\to\FF\VV ect_k$.

To some extent, the result might seem obvious. Thus the underlying functor of each homomorphism of rig categories from $\widehat{\FF\SSS et}$ to any other rig category $\SSS$ must preserve both the unit object and the sum, at least up to isomorphism. It follows that its action on objects is essentially given in a canonical way because every object in $\widehat{\FF\SSS et}$ is a finite sum of the unit object (the singleton). Less obvious, however, is the fact that the action on morphisms is also essentially determined by the axioms of a homomorphism of rig categories, so that there is an essentially unique way of mapping the symmetric group $S_n$, for each $n\geq 2$, into the group of automorphisms $Aut_\Ss(\underline{n})$ of the object $\underline{n}={\sf 1}+\stackrel{n)}{\cdots}+{\sf 1}$, with {\sf 1} the unit object of $\SSS$.~\footnote{Behind this uniqueness, there should be the coherence theorem for symmetric monoidal categories, although we made no use of it in the proof of the main theorem.} Moreover, there is also the point that, together with the underlying functor, giving a homomorphism from $\widehat{\FF\SSS et}$ to $\SSS$ also requires specifying the natural isomorphisms which take account of the preservation of $+$ and $\bullet$ up to isomorphisms. These natural isomorphisms must satisfy infinitely many coherence equations, and it is not  a priori clear that all possible choices for them actually define equivalent homomorphisms. 

Basic for the proof of the theorem will be working with a particular semistrict skeletal version of $\widehat{\FF\SSS et}$, as well as the fact that $\SSS$ can be assumed to be semistrict (i.e. such that many of the natural isomorphisms implicit in the structure of a rig category are identities). This allows us to describe quite explicitly the homomorphisms into $\SSS$, and to prove then that they are all indeed equivalent to a canonical one having the identity as its unique (2-)\,endomorphism.

Rig categories have been used since the late 1970s as sources of examples of $E_\infty$ ring spaces (see  Chapter~VI of \cite{May-1977}), and to define a sort of `2-K-theory' (see \cite{Baas-Dundas-Rognes-2004}). Our interest in rig categories is due to the fact that they constitute one of the basic inputs in the definition of the categorical analog of a (semi)module over a (semi)ring, the other one being a symmetric monoidal category to be acted on. More specifically, we are interested in the so-called $\FF\VV ect_k$-{\em module categories} as categorical analogs of the vector spaces. These are symmetric monoidal categories $\mathscr{M}=(\Mm,\oplus,{\sf 0})$ equipped with a categorical action of $\FF\VV ect_k$ on it. This is given by a homomorphism of rig categories from $\FF\VV ect_k$ to the rig category of endomorphisms of $\mathscr{M}$ (cf. Example~\ref{semianell_endomorfismes} below), and we would like to identify in more concrete terms the data that define such a $\FF\VV ect_k$-module category structure. Since $\widehat{\FF\SSS et}$ is biinitial in $\mathbf{RigCat}$, ebery $\mathscr{M}$ has a unique $\widehat{\FF\SSS et}$-module category structure up to equivalence, given by the essentially unique rig category homomorphism from $\widehat{\FF\SSS et}$ to the endomorphisms of $\mathscr{M}$ (cf. Example~\ref{morfisme_a_endomorfismes_M}). The point is that $\widehat{\FF\SSS et}$ embedds as a rig subcategory of $\FF in\VV ect_k$ through the free vector space construction. Hence it follows from our main theorem that part of any $\FF in\VV ect_k$-module category structure on $\mathscr{M}$ is canonically given. This should contribute to better understand the additional data that defines a $\FF in\VV ect_k$-module category structure on $\mathscr{M}$. 

\subsection{Outline of the paper and assumed background} In section~2 we review the definition of (symmetric) rig category, and the corresponding notions of 1- and 2-morphism making them the objects of a 2-category. Some examples are given, paying special attention to the symmetric rig category `canonically' associated to every distributive category. The section ends with the statement of the strictification theorem for (symmetric) rig categories. In Section~3 a detailed description is given of the semistrict version of $\widehat{\FF\SSS et}$ we shall use in this work. Finally, in Section~4 we prove that this 2-rig is biinitial in the 2-category of rig categories. Incidentally, we obtain explicit descriptions of the rig category homomorphisms from this semistrict version of $\widehat{\FF\SSS et}$ into any (semistrict) rig category $\SSS$. In a sense, these descriptions might have no interest because the category is contractible. However, as mentioned before, we are actually interested in the categories $\Hh om_{\mathbf{RigCat}}(\SSS',\SSS)$ for other rig categories $\SSS'$ containing $\widehat{\FF\SSS et}$ as a sub-2-rig, such as the rig category $\FF\VV ect_k$ of finite dimensional vector spaces over a given ground field. The point is that isomorphic objects in $\Hh om_{\mathbf{RigCat}}(\widehat{\FF\SSS et},\SSS)$ may no longer be equivalent when extended to homomorphisms from the whole rig category $\SSS'$. Thus having the above descriptions may be useful in the study of these categories of homomorphisms.

The reader is assumed to be familiar with the definitions of (symmetric) monoidal category and 2-category, and with the notions of (symmetric) monoidal functor and monoidal natural transformation between them. Good references for the basics of monoidal categories are, for instance, Chapters~1 and 3 of \cite{Aguiar-Mahajan-2010}, Chapter~XI of \cite{Kassel-1995}, or Chapter~2 of \cite{EGNO-2015}, and for the basics of 2-categories Chapter~7 of \cite{Borceux-1994-I}. The reader may also take a look to the standard text by MacLane \cite{MacLane-1998}.

\subsection{A few conventions about notation and terminology} Both sets and structured sets (such as monoids, rigs, etc) are denoted in the same way, namely by capital letters $A,B,C,\ldots$. Plain categories and functors between them are denoted by $\Aa,\Bb,\Cc,\ldots$, (symmetric) monoidal categories and (symmetric) monoidal functors between them by $\mathscr{A},\mathscr{B},\mathscr{C},\ldots$ and rig categories and rig category homomorphisms between them (both notions are defined below) by $\AAA,\BB,\CC,\ldots$. When we refer to concrete examples of categories, the same convention will be applied to the first letters. Thus $\Ss et$, $\mathscr{S}et$ and $\SSS et$ respectively denote the plain category of (small) sets, the category $\Ss et$ equipped with a monoidal structure, and the category $\Ss et$ equipped with a rig category structure. Finally, 2-categories are denoted by boldface letters $\Abf,\Bbf,\Cbf,\ldots$.

Given two monoidal categories $\mathscr{A}=(\Aa,\otimes,1,a,l,r)$ and $\mathscr{A}'=(\Aa',\otimes',1',a',l',r')$, and unless otherwise indicated, by a monoidal functor between them we mean a {\em strong} monoidal functor. The underlying functor of such a monoidal functor $\mathscr{F}:\mathscr{A}\to\mathscr{A}'$ is denoted by $F$, and the structural isomorphisms taking account of the preservation of the tensor product and unit objects up to natural isomorphism by $\varphi$ and $\varepsilon$, respectively. Thus $\mathscr{F}$ consists of the functor $F$ together with a natural isomorphism
\[
\varphi_{x,y}:F(x\otimes y)\stackrel{\cong}{\to} Fx\otimes' Fy
\]
in $\Aa'$ for each pair of objects $x,y\in\Aa$, and an isomorphism $\varepsilon:F1\to 1'$ also in $\Aa'$, all these isomorphisms satisfying the corresponding coherence axioms. 

By a {\em semistrict} symmetric monoidal category we shall mean what some authors call a strict symmetric monoidal category, i.e. a symmetric monoidal category whose associator and left and right unitors are identities, but not the commutator. The term strict will refer to the case when the commutator is also trivial. Every symmetric monoidal category is equivalent to a semistrict one, but not to a strict one (this is MacLane's coherence theorem for symmetric monoidal categories; cf. \cite{MacLane-1998}).

Finally, composition of (1-)morphisms (in a category or in a 2-category) is denoted by juxtaposition, and the identity of an object $x$ by $id_x$. In particular, for any objects $A,B,C$ and morphisms $f:A\to B$ and $g:B\to C$, the composite is denoted by $g\,f:A\to C$. Exceptionally, composition of functors is denoted with the symbol $\circ$\,.

\section{The 2-category of (symmetric) rig categories}
\label{categories_rig}

In this section we review the notion of (symmetric) rig category, motivating the required coherence axioms, and give various examples. Next appropriate notions of 1- and 2-cell are defined making (symmetric) rig categories the objects of a 2-category. In fact, there are various reasonable notions of 1-cell, and correspondingly various 2-categories of (symmetric) rig categories. We will stick to the {\em strong} notion of 1-cell (cf. Definition~\ref{morfisme_categories_semianell} below), although natural examples are given which are not strong. The section ends with the statement of the corresponding strictification theorem.

\subsection{Rig categories}
As recalled in the introduction, a rig (or semiring) is an abelian monoid $(S,0,+)$ equipped with an additional monoid structure $(\bullet\,,1)$ such that $\bullet$ distributes over $+$ from either side and $0\bullet x=x\bullet 0=0$ for each $x\in S$. Equivalently, $(\bullet\,,1)$ must be such that all left and right translation maps $L^x:y\mapsto x\bullet y$ and $R^x:y\mapsto y\bullet x$ are monoid endomorphisms of $(S,+,0)$. It follows that
\begin{itemize}
\item[(1)] $L^x\circ R^y=R^y\circ L^x$ and $L^x+R^y=R^y+L^x$ for each $x,y\in S$;
\item[(2)] $L^x\circ L^y=L^{x\bullet y}$ and $R^{y\bullet x}=R^x\circ R^y$ for each $x,y\in S$;
\item[(3)] $L^{x+y}=L^x+L^y$ and $R^{y+x}=R^y+R^x$ for each $x,y\in S$;
\item[(4)] $L^1=R^1={\sf 1}_S$;
\item[(5)] $L^0=R^0={\sf 0}_S$,
\end{itemize}
where ${\sf 1}_S$ and ${\sf 0}_S$ denote the identity and zero maps of $S$, respectively, and the sum of endomorphisms is pointwise defined. Notice that the first condition in (1), and each of the two conditions in (2) correspond to the associativity of $\bullet$\,. All of them are made explicit because they lead to different conditions in the categorified definition.

When the definition of a rig is categorified, the abelian monoid $(S,+,0)$ must be replaced by a symmetric monoidal category $\mathscr{S}=(\Ss,+,0,a,c,l,r)$, with $a,c,l,r$ the associativity, commutativity, and left and right unit natural isomorphisms, the monoid structure $(\bullet\,,1)$ by an additional monoidal structure $(\bullet\,,1,a',l',r')$ on $\Ss$, with $a',l',r'$ as before, and the distributivity and `absorbing' axioms by natural isomorphisms (for short, the symbol $\bullet$ between objects is omitted from now on)
\[
d_{x,y,z}:x(y+z)\stackrel{\cong}{\to}xy+xz,\quad d'_{x,y,z}:(x+y)z\stackrel{\cong}{\to}xz+yz
\]
\[
n_x:x0\stackrel{\simeq}{\rightarrow} 0,\quad n'_x:0x\stackrel{\simeq}{\rightarrow} 0
\]
making each of the left and right translation functors $L^x,R^x:\Ss\to\Ss$ into symmetric $+$-monoidal endofunctors $\mathscr{L}^x=(L^x,d_{x,-,-},n_x)$ and $\mathscr{R}^x=(R^x,d'_{-,-,x},n'_x)$ of $\mathscr{S}$, in the same way as the maps $L^x,R^x$ were monoid endomorphisms of $(S,+,0)$. In doing this, equalities (1)-(5) no longer hold strictly. We only have canonical natural isomorphisms between the corresponding pairs of functors, given by the natural isomorphisms
\begin{align*}
a'_{x,-,y}&:L^x\circ R^y\Rightarrow R^y\circ L^x,
\\
c_{x-,-y}&:L^x+R^y\Rightarrow R^y+L^x, 
\\
a'_{x,y,-}&:L^x\circ L^y\Rightarrow L^{xy},
\\
a'_{-,y,x}&:R^{yx}\Rightarrow R^x\circ R^y,
\\
d'_{x,y,-}&:L^{x+y}\Rightarrow L^x+L^y ,
\\
d_{-,y,x}&:R^{y+x}\Rightarrow R^y+R^x,
\\
l':L^1 &\,\Rightarrow {\bf 1}_{\Ss},\ 
r':R^1\Rightarrow {\bf 1}_\Ss,
\\
n:L^0 &\,\Rightarrow {\bf 0}_\Ss, \ n':R^0\Rightarrow {\bf 0}_\Ss,
\end{align*}
where ${\bf 1}_\Ss,{\bf 0}_\Ss$ respectively denote the identity and zero functors of $\Ss$. Moreover, the domain and codomain functors of all these isomorphisms are both symmetric $+$-monoidal. As a matter of fact, however, all these natural isomorphisms except $c_{x-,-y}$ need not be {\em monoidal}, and this must be explicitly required. Thus we are led to the following definition of a rig category.

\subsubsection{\sc Definition.}\label{categoria_rig}
A {\em rig category} is a symmetric monoidal category $\mathscr{S}=(\Ss,+,0,a,c,l,r)$ together with the following data:
\begin{itemize}
\item[(SC1)] an additional monoidal structure $(\bullet\,,1,a',l',r')$ (the {\em multiplicative} monoidal structure);
\item[(SC2)] two families of isomorphisms ({\em left and right distributors})
$$d_{x,y,z}:x(y+z)\stackrel{\simeq}{\rightarrow}x y+x z,$$
$$d'_{x,y,z}:(x+y) z\stackrel{\simeq}{\to}x z+y z$$
natural in $x,y,z\in\Ss$;
\item[(SC3)] two families of isomorphisms ({\em absorbing isomorphisms})
$$n_x:x0\stackrel{\simeq}{\rightarrow} 0,$$
$$n'_x:0x\stackrel{\simeq}{\rightarrow} 0$$
natural in $x\in\Ss$.
\end{itemize}
Moreover, these data must satisfy the following axioms:
\begin{enumerate}
\item[(SC4)] for each $x\in\Ss$ the triples $\mathscr{L}^x=(L^x,d_{x,-,-},n_x)$, $\mathscr{R}^x=(R^x,d'_{-,-,x},n'_x)$ are symmetric $+$-monoidal endofunctors of $\mathscr{S}$; more explicitly, this means that the diagrams
\begin{equation*}
\xymatrix{& x(y+(z+t))\ar[r]^{id_x\bullet a_{y,z,t}}\ar[d]_{d_{x,y,z+t}} & x((y+z)+t)\ar[d]^{d_{x,y+z,t}} \\ (A1.1) & xy+x(z+t)\ar[d]_{id_{xy}+d_{y,z,t}} & x(y+z)+xt\ar[d]^{d_{x,y,z}+id_{xt}} \\ & xy+(xz+xt)\ar[r]_{a_{xy,xz,xt}} & (xy+xz)+xt}\quad \xymatrix{(x+(y+z))t\ar[r]^{a_{x,y,z}\bullet id_t}\ar[d]_{d'_{x,y+z,t}} & ((x+y)+z)t\ar[d]^{d'_{x+y,z,t}} \\ xt+(y+z)t\ar[d]_{id_{xt}+d'_{y,z,t}} & (x+y)t+zt\ar[d]^{d'_{x,y,t}+id_{zt}} \\ xt+(yt+zt)\ar[r]_{a_{xt,yt,zt}} & (xt+yt)+zt}
\end{equation*}
\begin{equation*}
\xymatrix{(A1.2) & x(y+z)\ar[r]^{id_x\bullet c_{y,z}}\ar[d]_{d_{x,y,z}} & x(z+y)\ar[d]^{d_{x,z,y}} \\ & xy+xz\ar[r]_{c_{xy,xz}} & xz+xy}\quad \xymatrix{(y+z)x\ar[r]^{c_{y,z}\bullet id_x}\ar[d]_{d'_{y,z,x}} & (z+y)x\ar[d]^{d'_{z,y,x}} \\ yx+zx\ar[r]_{c_{yx,zx}} & zx+yx}
\end{equation*}
\begin{equation*}
\xymatrix{(A1.3) & x(0+y)\ar[r]^{d_{x,0,y}}\ar[d]_{id_x\bullet l_y} & x0+xy\ar[d]^{n_x+id_{xy}} \\ & xy & 0+xy\ar[l]^{l_{xy}}}\quad \xymatrix{(0+x)y\ar[r]^{d'_{0,x,y}}\ar[d]_{l_x\bullet id_y} & 0y+xy\ar[d]^{n'_y+id_{xy}} \\ xy & 0+xy\ar[l]^{l_{xy}}}
\end{equation*}
\begin{equation*}
\xymatrix{(A1.4) & x(y+0)\ar[r]^{d_{x,y,0}}\ar[d]_{id_x\bullet r_y} & xy+x0\ar[d]^{id_{xy}+n_x} \\ & xy & xy+0\ar[l]^{r_{xy}}}\quad \xymatrix{(x+0)y\ar[r]^{d'_{x,0,y}}\ar[d]_{r_x\bullet id_y} & xy+0y\ar[d]^{id_{xy}+n'_y} \\ xy & xy+0\ar[l]^{r_{xy}}}
\end{equation*}
commute for all objects $x,y,z,t\in\Ss$;

\item[(SC5)] for each $x,y\in\Ss$ the natural isomorphism $a'_{x,-,y}:L^x\circ R^y\Rightarrow R^y\circ L^x$ is a symmetric monoidal natural isomorphism $\mathscr{L}^x\circ\mathscr{R}^y\Rightarrow\mathscr{R}^y\circ\mathscr{L}^x$; more precisely, this means that the diagrams
\begin{equation*}
\xymatrix{(A2.1) & x((z+t)y)\ar[r]^{id_x\bullet d'_{z,t,y}}\ar[d]_{a'_{x,z+t,y}} & x(zy+ty)\ar[r]^{d_{x,zy,ty}} & x(zy)+x(ty)\ar[d]^{a'_{x,z,y}+a'_{x,t,y}}  \\ & (x(z+t))y\ar[r]_{d_{x,z,t}\bullet id_y} & (xz+xt)y\ar[r]_{d'_{xz,xt,y}} & (xz)y+(xt)y}
\end{equation*}
\begin{equation*}
\xymatrix{(A2.2) & x(0y)\ar[r]^{a'_{x,0,y}}\ar[d]_{id_x\bullet n'_y} & (x0)y\ar[r]^{n_x\bullet id_y} & 0y\ar[d]^{n'_y} \\ & x0\ar[rr]_{n_x} && 0}
\end{equation*}
commute for all objects $x,y,z,t\in\Ss$;

\item[(SC6)] for each $x,y\in\Ss$ the natural isomorphisms $a'_{x,y,-}:L^x\circ L^y\Rightarrow L^{xy}$ and $a'_{-,y,x}:R^{yx}\Rightarrow R^x\circ R^y$ are symmetric monoidal natural isomorphisms $\mathscr{L}^x\circ\mathscr{L}^y\Rightarrow\mathscr{L}^{xy}$ and $\mathscr{R}^{yx}\Rightarrow\mathscr{R}^x\circ\mathscr{R}^y$, respectively; more precisely, this means that the diagrams
\begin{equation*}
\xymatrix{(A3.1) & x(y(z+t))\ar[r]^{id_x\bullet d_{y,z,t}}\ar[d]_{a'_{x,y,z+t}} & x(yz+yt)\ar[r]^-{d_{x,yz,yt}} & x(yz)+x(yt)\ar[d]^{a'_{x,y,z}+a'_{x,y,t}} \\ & (xy)(z+t)\ar[rr]_{d_{xy,z,t}} && (xy)z+(xy)t }\quad \xymatrix{x(y0)\ar[r]^{a'_{x,y,0}}\ar[d]_{id_x\bullet n_{y}} & (xy)0\ar[d]^{n_{xy}} \\ x0\ar[r]_{n_x} & 0 }
\end{equation*}
\begin{equation*}
\xymatrix{(A3.2) & ((t+z)y)x\ar[r]^{d'_{t,z,y}\bullet id_x} & (ty+zy)x\ar[r]^-{d'_{ty,zy,x}} & (ty)x+(zy)x \\ & (t+z)(yx)\ar[u]^{a'_{t+z,y,x}}\ar[rr]_{d'_{t,z,yx}} && t(yx)+z(yx)\ar[u]_{a'_{t,y,x}+a'_{z,y,x}} } \quad \xymatrix{0(yx)\ar[r]^{a'_{0,y,x}}\ar[d]_{n'_{yx}} & (0y)x\ar[d]^{n'_y\bullet id_x} \\ 0 & 0x\ar[l]^{n'_x}}
\end{equation*}
commute for all $x,y,z,t\in\Ss$;

\item[(SC7)] for each $x,y\in\Ss$ the natural isomorphisms $d'_{x,y,-}:L^{x+y}\Rightarrow L^x+L^y$ and $d_{-,y,x}:R^{y+x}\Rightarrow R^y+R^x$ are symmetric monoidal natural isomorphisms $\mathscr{L}^{x+y}\Rightarrow\mathscr{L}^x+\mathscr{L}^y$ and $\mathscr{R}^{y+x}\Rightarrow\mathscr{R}^y+\mathscr{R}^x$, respectively; more precisely, this means that the diagrams
\begin{equation*}
\xymatrix{(A4.1) & (x+y)(z+t)\ar[r]^-{d'_{x,y,z+t}}\ar[d]_{d_{x+y,z,t}} & x(z+t)+y(z+t)\ar[r]^{d_{x,z,t}+d_{y,z,t}} & (xz+xt)+(yz+yt)\ar[d]^{v_{xz,xt,yz,yt}}  \\ & (x+y)z+(x+y)t\ar[rr]_{\ \ d'_{x,y,z}+d'_{x,y,t}} & & (xz+yz)+(xt+yt)}
\end{equation*}
\begin{equation*}
\xymatrix{(A4.2) & (x+y)0\ar[r]^{d'_{x,y,0}}\ar[d]_{n_{x+y}} & x0+y0\ar[d]^{n_x+n_y} \\ & 0 & 0+0\ar[l]^{l_0}}\quad \xymatrix{0(x+y)\ar[r]^{d_{0,x,y}}\ar[d]_{n'_{x+y}} & 0x+0y\ar[d]^{n'_x+n'_y} \\ 0 & 0+0\ar[l]^{l_0}}
\end{equation*}
commute for all objects $x,y,z,t\in\Ss$, where $v_{a,b,c,d}:(a+b)+(c+d)\stackrel{\cong}{\to}(a+c)+(b+d)$ in (A4.1) is the canonical isomorphism built from the associator $a$ and the symmetry $c$ of $\mathscr{S}$;

\item[(SC8)] the left and right unit isomorphisms $l':L^1\Rightarrow {\bf 1}_\Ss$ and $r':R^1\Rightarrow {\bf 1}_\Ss$ are symmetric monoidal natural isomorphisms $\mathscr{L}^1\Rightarrow {\bf 1}_\Ss$ and $\mathscr{R}^1\Rightarrow {\bf 1}_\Ss$, respectively; more explicitly, this means that the diagrams
\begin{equation*}
\xymatrix{(A5.1) & 1(x+y)\ar[rr]^{d_{1,x,y}}\ar[rd]_{l'_{x+y}}  & & 1 x+1 y\ar[ld]^{l'_x+l'_y} \\ & & x+y &} \quad \xymatrix{(x+y) 1\ar[rr]^{d'_{x,y,1}}\ar[rd]_{r'_{x+y}}  & & x 1+y 1\ar[ld]^{r'_x+l'_y} \\ & x+y &}
\end{equation*}
commute for all objects $x,y\in\Ss$, and the following equalities hold:
\begin{equation*}
(A5.2)\qquad n_1=l'_0,\quad n'_1=r'_0;
\end{equation*}

\item[(SC9)] the left and right null isomorphisms $n':L^0\Rightarrow {\bf 0}_\Ss$ and $n:R^0\Rightarrow {\bf 0}_\Ss$ are symmetric monoidal natural isomorphisms $\mathscr{L}^0\Rightarrow {\bf 0}_\Ss$ and $\mathscr{R}^0\Rightarrow {\bf 0}_\Ss$, respectively; more explicitly, this means that
$$n_0=n'_0.$$
\end{enumerate}

\noindent
For short, we shall write $\SSS$ to denote the whole data $(\Ss,+,\bullet\,,0,1,a,c,l,r,a',l',r',d,d',n,n')$ defining a rig category. The objects $0$ and $1$ will be respectively called the {\em zero} and {\em unit} objects of $\SSS$. $\SSS$ will be called a {\em 2-rig} when the underlying category $\Ss$ is a groupoid.

\medskip
When convenient, and in order to distinguish the additive and the multiplicative monoidal structures on a rig category we shall write $\mathscr{S}^+=(\Ss,+,0,a,c,l,r)$, and $\mathscr{S}^{\,\bullet}=(\Ss,\bullet\,,1,a',l',r')$.

\subsubsection{\sc Definition.}
A rig category $\SSS$ is called {\em left} (resp. {\em right}) {\em semistrict} when all structural isomorphisms except $c$ and the right distributor $d'$ (resp. the left distributor $d$) are identities. It is called {\em semistrict} when it is either left or right semistrict, and {\em strict} when $c$ and both distributors are also identities.

\subsubsection{\sc Remark.}
Rigs can be more compactly defined as the one-object categories enriched over the symmetric monoidal closed category of abelian monoids with the usual tensor product of abelian monoids. Similarly, rig categories should correspond to one-object categories enriched over a suitable symmetric monoidal closed category of symmetric monoidal categories. Such enriched categories are considered by Guillou \cite{Guillou-2010} under the name of one-object {\em {\bf SMC}-categories}, although he avoids describing explicitly the symmetric monoidal structure on the category {\bf SMC} of symmetric monoidal categories.

\subsection{Symmetric rig categories.}\label{categories_rig_simetriques}
A rig $(S,+,\bullet,0,1)$ is commutative when the monoid $(S,\bullet,1)$ is abelian, or equivalently when $L^x=R^x$ for each $x\in S$. When this condition is categorified, the abelian monoid structure $(\bullet,1)$ becomes a symmetric monoidal structure $(\bullet,1,a',c',l',r')$ on $\Ss$, and instead of the equality $L^x=R^x$ we now have the natural isomorphism $c'_{x,-}:L^x\Rightarrow R^x$. Once more, this isomorphism may not be a symmetric monoidal natural isomorphism $\mathscr{L}^x\Rightarrow\mathscr{R}^x$, and this must be required explicitly. This leads us to the following definition.

\subsubsection{\sc Definition.}\label{categoria_rig_simetrica}
A {\em symmetric rig category} is a rig category $\SSS$ together with a family of natural isomorphisms $c'_{x,y}:x y\stackrel{\simeq}{\to}y x$, for each $x,y\in\Ss$, called the {\em multiplicative commutators}, such that
\begin{itemize}
\item[(CSC1)] $(\Ss,\bullet,1,a',c',l',r')$ is a symmetric monoidal category,
\item[(CSC2)] for each $x\in\Ss$, the natural isomorphism $c'_{x,-}:L^x\Rightarrow R^x$ is a symmetric monoidal natural isomorphism $\mathscr{L}^x\Rightarrow\mathscr{R}^x$; more precisely, this means that the diagrams
\[
\xymatrix{(x+y) z\ar[r]^{d'_{z;x,y}}\ar[d]_{c'_{x+y,z}} & x z+y z\ar[d]^{c'_{x,z}+c'_{y,z}} \\ z(x+y)\ar[r]_{d_{z;x,y}} & z x+z y}\qquad 
\xymatrix{x0\ar[rr]^{c'_{x,0}}\ar[dr]_{n_x} && 0x\ar[ld]^{n'_x} \\ & 0 &}
\]
commute for all objects $x,y,z\in\Ss$.
\end{itemize}
A symmetric rig category is called {\em semistrict} (resp. {\em strict}) when the underlying rig category $\SSS$ is semistrict (resp. strict and with $c'$ trivial). ~\footnote{\,Some people, mostly those working on the K-theory of this type of categories, call the semistrict symmetric rig categories {\em bipermutative categories} because a semistrict symmetric monoidal category is also called a {\em permutative category}.} It is called a {\em symmetric 2-rig} when the underlying rig category is a 2-rig.

\subsubsection{\sc Remark.}
The previous definition coincides with the structure described by Laplaza in \cite{Laplaza-1972} except that in Laplaza's paper the distributors are only required to be monomorphisms. The correspondence between Laplaza's axioms and the axioms in Definitions~\ref{categoria_rig} and \ref{categoria_rig_simetrica} goes as follows: 
\begin{itemize}
\item[(i)] (A1.1)-(A1.4) correspond to Laplaza's axioms (I), (III)-(V) and (XIX)-(XXII), 
\item[(ii)] (A2.1)-(A2.2) correspond to Laplaza's axioms (VIII) and (XVII),
\item[(iii)] (A3.1)-(A3.2) correspond to Laplaza's axioms (VI)-(VII), (XVI) and (XVIII),
\item[(iv)] (A4.1)-(A4.2) correspond to Laplaza's axioms (IX) and (XI)-(XII),
\item[(v)] (A5.1)-(A5.2) correspond to Laplaza's axioms (XIII)-(XIV) and (XXIII)-(XXIV),
\item[(vi)] (SC9) corresponds to Laplaza's axiom (X), and
\item[(vii)] (CSC2) correspond to Laplaza's axioms (II) and (XV).
\end{itemize}

\subsection{Some examples of (symmetric) rig categories}
Describing a rig category requires specifying the data $\Ss,+,\bullet,0,1,a,c,l,r,a',l',r',d,d',n,n'$, and checking that they satisfy all of the above axioms. Usually, this is long and tedious. Hence in this subsection we just mention a few standard examples of rig categories without getting into the details. The particular type of rig categories we are interested in is discussed in more detail in \S~\ref{categories_distributives}.

\subsubsection{\sc Example.}
Every rig $S=(S,+,\bullet,0,1)$ can be thought of as a 2-rig with only identity morphisms, and all required isomorphisms trivial. They are symmetric 2-rigs when $S$ is commutative.

\subsubsection{\sc Example.}\label{core}
If $\SSS$ is a rig category, and $\widehat{\Ss}$ is the groupoid with the same objects as $\Ss$ and only the isomorphisms between them as morphisms then $\widehat{\Ss}$ inherits by restriction a canonical rig category structure. The 2-rig so obtained is denoted by $\widehat{\SSS}$. It is a symmetric 2-rig when $\SSS$ is a symmetric rig category, and (left,right) semistrict when $\SSS$ is so.

\subsubsection{\sc Example.}\label{categories_semianell_no_cartesianes}
Every symmetric monoidal closed category (in particular, every cartesian closed category) with finite coproducts is canonically a symmetric rig category with $+$ and $\bullet$ respectively given by the categorical coproduct $\sqcup$ and the tensor product $\otimes$. Closedness is necessary to ensure that the tensor product indeed distributes over coproducts. Thus denoting the internal homs by $[-,-]$ we have
\begin{align*}
Hom((x\sqcup y)\otimes z,t)&\cong Hom(x\sqcup y,[z,t])
\\ &\cong Hom(x,[z,t])\otimes Hom(y,[z,t])
\\ &\cong Hom(x\otimes z,t)\otimes Hom(y\otimes z,t)
\\ &\cong Hom((x\otimes z)\sqcup(y\otimes z),t)
\end{align*}
for every objects $x,y,z,t$. Then the existence of the isomorphism $d'_{x,y,z}$ follows from the Yoneda lemma. A similar argument gives the isomorphism $d_{x,y,z}$. Symmetric rig categories of this type include those associated to the three cartesian closed categories $\Ss et$ of sets and maps (and its full subcategory $\Ff\Ss et$ with objects the finite sets), $\Ss et_G$ of $G$-sets and homomorphisms of $G$-sets for any group $G$ (and its full subcategory $\Ff\Ss et_{G}$ with objects the finite $G$-sets), and $\Cc at$ of (small) categories and functors, and those associated to the two non-cartesian symmetric monoidal closed categories $\Vv ect_k$ of vector spaces over a field $k$ and $k$-linear maps (and its full subcategory $\Ff\Vv ect_k$ with objects the finite dimensional vector spaces), and $\Rr ep_k(G)$ of $k$-linear representations of a group $G$ and homomorphisms of representations for any group $G$ (and its full subcategory $\Ff\Rr ep_k(G)$ with objects the finite dimensional representations).

\subsubsection{\sc Example.} \label{semianell_endomorfismes}
The set of endomorphisms of every abelian monoid is canonically a rig, with the sum of monoid endomorphisms defined pointwise, and the product given by the composition. Similarly, for every symmetric monoidal category $\mathscr{M}=(\Mm,\oplus,0_\Mm,\asf,\csf,\lsf,\rsf)$ the category $\Ee nd\,(\mathscr{M})$ of symmetric monoidal endofunctors of $\mathscr{M}$, and symmetric monoidal  natural transformations between them is canonically a (non-symmetric) rig category, with the additive symmetric monoidal structure given by the pointwise sum of symmetric monoidal endofunctors and monoidal natural transformations, and with the multiplicative monoidal structure given by their composition. The zero object is the zero functor ${\bf 0}_\Mm$, and the canonical isomorphisms $a,c,l,r$ are pointwise given by the corresponding isomorphisms in $\mathscr{M}$. In particular, $\Ee nd(\mathscr{M})$ is semistrict (strict) symmetric monoidal category when $\mathscr{M}$ is so. The multiplicative monoidal structure is always strict, with the identity functor of $\Mm$ as unit object. Moreover, if we define $\mathscr{F}\bullet\mathscr{F}'=\mathscr{F}'\circ\mathscr{F}$, the left distributors are all trivial, while the right distributors
\[
d'_{\mathscr{F}_1,\mathscr{F}_2,\mathscr{F}_3}:\mathscr{F}_3\circ(\mathscr{F}_1+\mathscr{F}_2)\Rightarrow\mathscr{F}_3\circ\mathscr{F}_1+\mathscr{F}_3\circ\mathscr{F}_2
\]
are given by the monoidal structure of $\mathscr{F}_3$. As to the absorbing isomorphisms $n_{\mathscr{F}},n'_{\mathscr{F}}$, they are all trivial. In particular, the rig category so defined $\EE nd(\mathscr{M})$ is left semistrict when $\mathscr{M}$ is semistrict, and a 2-rig when $\Mm$ is a groupoid.

\subsection{Distributive categories}\label{categories_distributives}
Recall that a {\em distributive category} is a cartesian and cocartesian category such that the canonical map $xy+xz\to x(y+z)$ is invertible for each objects $x,y,z$. Distributive categories generalize the cartesian closed categories with finite coproducts of Example~\ref{categories_semianell_no_cartesianes}, and as these have a ``canonical'' symmetric rig category structure. In this paragraph this structure is described in detail.

Let us first recall that every cocartesian category $\Cc$ (i.e. a category $\Cc$ with all finite coproducts) has a ``canonical'' symmetric monoidal structure associated to the choice of a particular coproduct $(x+y,\iota^1_{x,y},\iota^2_{x,y})$ for each ordered pair of objects $(x,y)$, and a particular initial object $0$. It is given as follows:
\begin{itemize}
\item[(D1)] the tensor product $\Cc\times\Cc\to\Cc$ is given on objects $(x,y)$ and morphisms $(f,g):(x,y)\to(x',y')$ by
\[
(x,y)\mapsto x+y,\quad (f,g)\mapsto f+g,
\]
where $f+g$ is the morphism uniquely determined by the diagram
\[
\xymatrix{
x\ar[r]^-{\iota^1_{x,y}}\ar[d]_-f & x+y\ar@{.>}[d]^-{\exists !\,f+g} & y\ar[l]_-{\iota^2_{x,y}}\ar[d]^-g
\\
x'\ar[r]_-{\iota^1_{x',y'}} & x'+y' &  y'\ar[l]^-{\iota^2_{x',y'}}
}
\]
and the universal property of $x+y$;

\item[(D2)] the unit object is the chosen initial object $0$; 
 
\item[(D3)] for every objects $x,y,z\in\Cc$ the associator $a_{x,y,z}$ is the morphism uniquely determined by the left hand side diagram
\[
\xymatrix{
x\ar[d]_-{\iota^1_{x,y}}\ar[r]^-{\iota^1_{x,y+ z}} 
& 
x+(y+z)\ar@{.>}[d]_-{a_{x,y,z}}
&
y+z\ar[l]_-{\ \ \iota^2_{x,y+ z}}\ar[ld]^-{\iota^2_{x,y}\,+\, id_z}
\\
x+y\ar[r]_-{\iota^1_{x+y,z}} & (x+y)+z  & 
}
\qquad
\xymatrix{
x+y\ar[r]^-{\iota^1_{x+ y,z}}\ar[rd]_-{id_x\,+\,\iota^1_{y,z}} 
& 
(x+y)+z\ar@{.>}[d]^-{a^{-1}_{x,y,z}}
&
z\ar[l]_-{\iota^2_{x+y,z}}\ar[d]^-{\iota^2_{y,z}}
\\
& x+(y+z) & y+ z\ar[l]^-{\iota^2_{x,y+ z}}
}
\]
and the universal properties of $x+(y+z)$ (it is indeed invertible with inverse the morphism uniquely determined by the right hand side diagram and the universal property of $(x+y)+z$);

\item[(D4)] for every objects $x,y\in\Cc$ the commutator $c_{x,y}$ is the morphism uniquely determined by the diagram
\[
\xymatrix{
x\ar[r]^-{\iota^1_{x,y}}\ar[dr]_-{\iota^2_{y,x}} &
x+y\ar@{.>}[d]^-{c_{x,y}} & 
y\ar[l]_-{\iota^2_{x,y}}\ar[ld]^{\iota^1_{y,x}} 
\\
& y+x & 
}
\]
and the universal property of $x+y$ (it is indeed invertible with inverse $c^{-1}_{x,y}=c_{y,x}$);

\item[(D5)] for every object $x\in\Cc$ the left and right unitors $l_x,r_x$ are the morphisms uniquely determined by the diagrams
\[
\xymatrix{
0\ar[r]^-{!}\ar[dr]_-{!} &
0+x\ar@{.>}[d]_-{l_x} & 
x\ar[l]_-{\iota^2_{0,x}}\ar[ld]^{id_x} 
\\
& x & 
}
\qquad
\xymatrix{
x\ar[r]^-{\iota^1_{x,0}}\ar[dr]_-{id_x} &
x+0\ar@{.>}[d]^-{r_x} & 
0\ar[l]_-{!}\ar[ld]^{!} 
\\
& x & 
}
\]
and the universal properties of $0+x$ and $x+0$ (they are indeed invertible with inverses $l^{-1}_x,r^{-1}_x$ the morphisms $\iota^2_{0,x},\iota^1_{x,0}$ respectively).
\end{itemize}
Different choices of binary coproducts and initial object lead to different but equivalent symmetric monoidal structures on $\Cc$ and hence, we may indeed speak of the ``canonical'' symmetric monoidal structure on each cocartesian category $\Cc$. Similarly, every cartesian category $\Cc$ (i.e. a category $\Cc$ with all finite products) is ``canonically'' a symmetric monoidal category with the associator $a'$, commutator $c'$, and left and right unitors $l',r'$ defined by the dual diagrams for some particular choices of binary products and final object.

When $\Cc$ is both cartesian and cocartesian, these two symmetric monoidal structures $(+,0,a,c,l,r)$ and $(\bullet,1,a',c',l',r')$ are related by the natural {\em left} and {\em right distributor maps}
\begin{align*}
\overline{d}_{x,y,z}&:xy+xz\to x(y+z) \\ \overline{d}\,'_{x,y,z}&:xz+yz\to (x+y)z
\end{align*}
uniquely determined by the diagrams
\[
\xymatrix{
xy\ar[r]^-{\iota^1_{xy,xz}}\ar[dr]_-{id_x\bullet\,\iota^1_{y,z}} &
xy+xz\ar@{.>}[d]^-{\overline{d}_{x,y,z}} & 
xz\ar[l]_-{\iota^2_{xy,xz}}\ar[ld]^{id_x\bullet\,\iota^2_{y,z}} 
\\
& x(y+z) & 
}
\quad
\xymatrix{
xz\ar[r]^-{\iota^1_{xz,yz}}\ar[dr]_-{\iota^1_{x,y}\bullet\,id_z} &
xz+yz\ar@{.>}[d]^-{\overline{d}\,'_{x,y,z}} & 
yz\ar[l]_-{\iota^2_{xz,yz}}\ar[ld]^{\iota^2_{x,y}\bullet\,id_z} 
\\
& (x+y)z & 
}
\]
and the universal properties of the coproducts $xy+xz$ and $xz+yz$. In fact, these distributors are not independent of each other. Instead, it may be shown that they are related bu the commutativity of the diagram
\[
\xymatrix{
xz+yz\ar[r]^-{\overline{d}\,'_{x,y,z}}\ar[d]_-{c'_{x,z}+c'_{y,z}} & (x+y)z\ar[d]^-{c'_{x+y,z}}
\\
zx+zy\ar[r]_-{\overline{d}_{z,x,y}} & z(x+y)\,.
}
\]
In particular, $\overline{d}\,'_{x,y,z}$ is invertible if and only if $\overline{d}_{z,x,y}$ is invertible.

The distributors $\overline{d}_{x,y,z},\overline{d}\,'_{x,y,z}$ are in general non-invertible, and as said before, $\Cc$ is called {\em distributive} precisely when $\overline{d}_{x,y,z}$ (or equivalently, $\overline{d}\,'_{x,y,z}$) is invertible for every objects $x,y,z$. In fact, for a cartesian and cocartesian category to be distributive it is enough that there exists {\em any} natural isomorphism $xy+xz\cong x(y+z)$ (see \cite{Lack-2012}). The point is that, when $\Cc$ is distributive, there are also canonical isomorphisms $x0\cong 0\cong 0x$ for each $x\in\Cc$ such that the whole structure makes $\Cc$ into a symmetric rig category. These are the maps $\pi^1_{0,x}:0x\to 0$ and $\pi^2_{x,0}:x0\to 0$, whose respective inverses are the unique maps $0\to 0x$ and $0\to x0$. Thus we have the following well known result (at least, among experts)

\subsubsection{\sc Proposition.}\label{estructura_categoria_rig_simetrica}
Every distributive category $\Cc$ equipped with the above symmetric monoidal structures $(+,0,a,c,l,r)$ and $(\bullet,1,a',c',l',r')$, and with the distributors and absorbing isomorphisms given by
\[
d_{x,y,z}=(\overline{d}_{x,y,z})^{-1},\quad d'_{x,y,z}=(\overline{d}\,'_{x,y,z})^{-1},
\quad n_x=\pi^2_{x,0},\quad n'_x=\pi^1_{0,x}
\]
is a symmetric rig category. 

\medskip
A standard example of a distributive category is the category $\Ff\Ss et$ of finite sets and maps between them. In Section~3 the symmetric rig category structure of a skeleton of it corresponding to a particular choice of products and coproducts is explicitly described.

\subsection{The 2-category of (symmetric) rig categories}

(Symmetric) rig categories are the objects of a 2-category. In fact, there are various useful notions of 1-cell between (symmetric) rig categories, associated to the various notions of 1-cell between (symmetric) monoidal categories, either lax, colax, bilax, strong or strict (symmetric) monoidal functor (see \cite{Aguiar-Mahajan-2010}). Moreover, we may also consider 1-cells whose character is different for the additive and the multiplicative monoidal structures. For instance, a 1-cell may be additively strong and multiplicatively colax, and examples of these mixed kind actually arise in some natural situations. Thus there are actually various 2-categories of (symmetric) rig categories. Although we shall define the various types of 1-cell, and give examples of various types, at the end we shall restrict to the strong morphisms and the associated 2-categories.

Recall that, given two rigs $S=(S,+,\bullet,0,1)$ and $\tilde{S}=(\tilde{S},\tilde{+},\tilde{\bullet},\tilde{0},\tilde{1})$, a rig homomorphism from $S$ to $\tilde{S}$ is a map $f:S\to\tilde{S}$ such that $f$ is both a monoid homomorphism from $(S,+,0)$ to $(\tilde{S},\tilde{+},\tilde{0})$, and a monoid homomorphism from $(S,\bullet,1)$ to $(\tilde{S},\tilde{\bullet},\tilde{1})$. It follows that $f$ is such that
\begin{align}
f\circ L^x&=\tilde{L}^{fx}\circ f \label{eq1} \\ f\circ R^x&=\tilde{R}^{fx}\circ f \label{eq2}
\end{align}
for each $x\in S$, where $L^x,R^x$ and $\tilde{L}^{fx},\tilde{R}^{fx}$ respectively denote the left and right translation maps of $S$ and $\tilde{S}$. In categorifying this definition, the map $f$ must be replaced by a functor $F:\Ss\to\tilde{\Ss}$ together with a pair $(\varphi^+,\varepsilon^+)$ making it a symmetric monoidal functor $\mathscr{S}^+\to\tilde{\mathscr{S}}^+$ of some type, and a pair $(\varphi^{\,\bullet},\varepsilon^{\,\bullet})$ making it a monoidal functor $\mathscr{S}^{\,\bullet}\to\tilde{\mathscr{S}}^{\,\bullet}$ of perhaps a different type. When done, equalities (\ref{eq1})-(\ref{eq2}) no longer hold. Instead, we just have natural transformations between the involved functors. For instance, in case the multiplicative monoidal structure is lax, we have the natural morphisms
\begin{align*}
\varphi^{\,\bullet}_{x,-}&:\tilde{L}^{Fx}\circ F\Rightarrow F\circ L^x\,, \\ \varphi^{\,\bullet}_{-,x}&:\tilde{R}^{Fx}\circ F\Rightarrow F\circ R^x.
\end{align*}
Moreover, the domain and codomain functors of both transformations are always symmetric monoidal of some kind, with monoidal structure given by the monoidal structure on $F$ and the distributors. However, $\varphi^{\,\bullet}_{x,-}$ and $\varphi^{\,\bullet}_{-,x}$ need not be monoidal. Thus we are naturally led to the following notions of 1-cell between rig categories.

\subsubsection{\sc Definition.}\label{morfisme_categories_semianell}
Let be given two rig categories $\SSS$ and $\tilde{\SSS}$. A {\em (colax,lax) morphism of rig categories} from $\SSS$ to $\tilde{\SSS}$ is a functor $F:\Ss\to\tilde{\Ss}$ together with the following data:
\begin{itemize}
\item[(HSC1)] an additive symmetric colax monoidal structure $(\varphi^+,\varepsilon^+)$ on $F$, and
\item[(HSC2)] a multiplicative lax monoidal structure $(\varphi^{\,\bullet},\varepsilon^{\,\bullet})$ on $F$.
\end{itemize}
Moreover, these data must satisfy the following axioms:
\begin{itemize}
\item[(HSC3)] for each $x\in\Ss$ the natural morphism $\varphi^{\,\bullet}_{x,-}:\tilde{L}^{Fx}\circ F\Rightarrow F\circ L^x$ is a symmetric monoidal natural isomorphism $\tilde{\mathscr{L}}^{Fx}\circ\mathscr{F}^a\Rightarrow\mathscr{F}^a\circ\mathscr{L}^x$, i.e. the diagrams
\[
\xymatrix{
F(x(y+z))\ar[r]^-{F(d_{x,y,z})} & F(xy+xz)\ar[r]^-{\varphi^+_{xy,xz}} & F(xy)\,+\,F(xz)
\\
Fx\,F(y+z)\ar[u]^{\varphi^{\,\bullet}_{x,y+z}}\ar[r]_-{id_{Fx}\bullet\,\varphi^+_{y,z}} & Fx\,(Fy\,+\,Fz)\ar[r]_-{\tilde{d}_{Fx,Fy,Fz}}  & Fy\,Fx\,+\,Fz\,Fx\ar[u]_{\varphi^{\,\bullet}_{x,y}\,+\,\varphi^{\,\bullet}_{x,z}} 
} \quad
\xymatrix{
F(x0)\ar[r]^-{F(n_x)} & F0\ar[r]^-{\varepsilon^+} & \tilde{0}
\\
Fx\,F0\ar[u]^{\varphi^{\,\bullet}_{x,0}}\ar[r]_-{id_{Fx}\bullet\,\varepsilon^+} & (Fx)\,\tilde{0}\ar[r]_-{\tilde{n}_{Fx}} & \tilde{0}\ar@{=}[u]
}
\]
commute for each objects $y,z\in\Ss$;
\item[(HSC4)] for each $x\in\Ss$ the natural isomorphism $\varphi^{\,\bullet}_{-,x}:\tilde{R}^{Fx}\circ F\Rightarrow F\circ R^x$ is a symmetric monoidal natural isomorphism $\tilde{\mathscr{R}}^{Fx}\circ\mathscr{F}^a\Rightarrow\mathscr{F}^a\circ\mathscr{R}^x$, i.e. the diagrams
\[
\xymatrix{F((y+z)x)\ar[r]^-{F(d'_{y,z,x})} & F(yx+zx)\ar[r]^-{\varphi^+_{yx,zx}} & F(yx)\,+\,F(zx)
\\
F(y+z)\,Fx\ar[u]^{\varphi^{\,\bullet}_{y+z,x}}\ar[r]_-{\varphi^+_{y,z}\bullet\, id_{Fx}} & (Fy\,+\,Fz)\,Fx\ar[r]_-{\tilde{d}'_{Fy,Fz,Fx}} & Fy\,Fx\,+\,Fz\,Fx\ar[u]_{\varphi^{\,\bullet}_{y,x}\,+\,\varphi^{\,\bullet}_{z,x}}
} \quad
\xymatrix{F(0x)\ar[r]^-{F(n'_x)} & F0\ar[r]^-{\varepsilon^+} & \tilde{0}
\\
F0\,Fx\ar[u]^{\varphi^{\,\bullet}_{0,x}}\ar[r]_-{\varepsilon^+\bullet\,id_{Fx}} & \tilde{0}\,(Fx)\ar[r]_-{\tilde{n}'_{Fx}} & \tilde{0}\ar@{=}[u]
}
\]
commute for each objects $y,z\in\Ss$.
\end{itemize}

\noindent
For any other choices of $\alpha,\beta\in\{lax,\,colax,\,strong,\,strict\}$, $(\alpha,\beta)$ {\em morphisms of rig categories} are defined similarly. When $\alpha=\beta$, we shall speak of an $\alpha$ morphism. Finally, when $\SSS,\tilde{\SSS}$ are symmetric rig categories, an $(\alpha,\beta)$ {\em morphism of symmetric rig categories} from $\SSS$ to $\tilde{\SSS}$ is an $(\alpha,\beta)$-morphism such that $(\varphi^{\,\bullet},\varepsilon^{\,\bullet})$ is a symmetric $\beta$-monoidal structure on $F$, and an $\alpha$ {\em morphism of symmetric rig categories} is an $\alpha$-morphism such that $(\varphi^{\,\bullet},\varepsilon^{\,\bullet})$ is a symmetric $\alpha$-monoidal structure on $F$.

\medskip
For short, we shall denote by $\FF$ the whole data $(F,\varphi^+,\varepsilon^+,\varphi^{\,\bullet},\varepsilon^{\,\bullet})$ defining a morphism of rig categories of any type.

\subsubsection{\sc Remark.}
A strong morphism of rig categories corresponds to the one-object case of the more general notion of a {\em {\bf SMC}-functor} between {\bf SMC}-categories introduced by Guillou (\cite{Guillou-2010}, Definition~4.1).

\subsubsection{\sc Example.}
Every rig homomorphism between two rigs $S$ and $\tilde{S}$ is a strict morphism between the associated discrete 2-rigs, and conversely.

\subsubsection{\sc Example.}
For every rig category $\SSS$ the inclusion functor $J:\widehat{\Ss}\hookrightarrow\Ss$ of the underlying groupoid $\widehat{\Ss}$ into $\Ss$ is a strict morphism of rig categories $\JJ:\widehat{\SSS}\hookrightarrow\SSS$ (cf. Example~\ref{core}).

\subsubsection{\sc Example.}
Let $\SSS et_G$ be the symmetric rig category of $G$-sets for some group $G$ (cf. Example~\ref{categories_semianell_no_cartesianes}). Then the forgetful functor $U_G:\Ss et_G\to\Ss et$ is a strict morphism of symmetric rig categories. However, its left adjoint $J_G:\Ss et\to\Ss et_G$, mapping each set $X$ to $X\times G$ with $G$-action given by $g'(x,g)=(x,g'g)$, and each map $f:X\to Y$ to the morphism of $G$-sets $f\times id_G:X\times G\to Y\times G$, is canonically just a strong-colax morphism $\JJ_G:\SSS et\to\SSS et_G$. The additive strong monoidal structure is given by the canonical right distributors
\[
\varphi^+_{X,Y}=d'_{X,Y,G}:(X\sqcup Y)\times G\stackrel{\cong}{\to} (X\times G)\sqcup(Y\times G)\,,
\]
together with the unique map $\varepsilon^+:\emptyset\times G\stackrel{\cong}{\to} \emptyset$, while the multiplicative colax structure is given by the canonical non-invertible morphisms of $G$-sets
\[
\varphi^{\,\bullet}_{X,Y}:(X\times Y)\times G\to (X\times G)\times(Y\times G)
\]
defined by $((x,y),g)\mapsto((x,g),(y,g))$, together with the unique map $\varepsilon^{\,\bullet}:\{*\}\times G\to\{*\}$.

\subsubsection{\sc Example.}
Let $\VV ect_k$ be the symmetric rig category of vector spaces over a given field $k$ (cf. Example~\ref{categories_semianell_no_cartesianes}). Then the forgetful functor $U_k:\Vv ect_k\to\Ss et$ is a lax morphism of symmetric rig categories $\UU_k:\VV ect_k\to\SSS et$ with additive lax monoidal structure given by the canonical maps $\varphi^+_{V,W}:V\sqcup W\to V\times W$ defined by $v\mapsto(v,0)$ and $w\mapsto (0,w)$, together with the canonical map $\varepsilon^+:\emptyset\to \{0\}$, while the multiplicative lax structure is given by the canonical maps $\varphi^{\,\bullet}_{V,W}:V\times W\to V\otimes_k W$ given by $(v,w)\mapsto v\otimes w$, and the map $\varepsilon^\times:\{*\}\to k$ sending $*$ to the unit $1\in k$. By constrast, its left adjoint $J_k:\Ss et\to\Vv ect_k$, mapping each set $X$ to the vector space $k[X]$ spanned by $X$, is canonically a strong morphism of symmetric rig categories $\JJ_k:\SSS et\to\VV ect_k$ with $\varphi^+_{X,Y}$, $\varepsilon^+$, $\varphi^{\,\bullet}_{X,Y}$, $\varepsilon^{\,\bullet}$ the usual natural isomorphisms $k[X\sqcup Y]\cong k[X]\oplus k[Y]$, $k[\emptyset]\cong 0$, $k[X\times Y]\cong k[X]\otimes k[Y]$, and $k[\{*\}]\cong k$.

\medskip
From now on,  we shall restrict to strong morphisms of (symmetric) rig categories, and they will be called {\em homomorphisms}.

\subsubsection{\sc Definition.} \label{transformacio_rig}
Let $\SSS,\tilde{\SSS}$ be two (symmetric) rig categories, and $\FF_1,\FF_2:\SSS\to\tilde{\SSS}$ two homomorphisms between them. A {\em rig transformation} from $\FF_1$ to $\FF_2$ is a natural transformation $\xi:F_1\Rightarrow F_2$ that is both $+$-monoidal and $\bullet$\,-monoidal.

\subsubsection{\sc Remark.} 
There is a more general notion of 2-cell $\FF_1\Rightarrow\FF_2$ corresponding to Guillou's definition of monoidal transformation between {\bf SMC}-functors whose domain and codomain {\bf SMC}-categories have only one object and hence, are rig categories (\cite{Guillou-2010}, Definition~4.2). It consists of an object $\tilde{x}$ in $\tilde{\Ss}$ together with a family of natural morphisms $\eta_x:(F_2\,x)\,\tilde{x}\to\tilde{x}\,(F_1x)$ in $\tilde{\Ss}$, labelled by the objects $x$ in $\Ss$, satisfying appropriate conditions. This is analogous to the existence of a more general notion of 2-cell between (symmetric) monoidal functors, corresponding to the pseudonatural transformations between them when viewed as pseudofunctors between one-object bicategories. Then the previous notion of rig transformation is to be thought of as the analog in the rig category setting of Lack's icons \cite{Lack-2010}.

\medskip
Rig categories together with the rig category homomorphisms as 1-cells, and the rig transformations between these as 2-cells constitute a 2-category $\mathbf{RigCat}$. The various compositions of 1- and 2-cells are defined in the obvious way. Similarly, symmetric rig categories with the symmetric rig category homomorphisms, and the rig transformations as 2-cells also constitute a 2-category $\mathbf{SRigCat}$. Notice that, unlike the category $\Cc\Rr ig$ of commutative rigs, which is a full subcategory of $\Rr ig$, $\mathbf{SRigCat}$ is not a full sub-2-category of $\mathbf{RigCat}$ because being a {\em symmetric} rig category is not a property-like structure. A given rig category can be symmetric in various non-equivalent ways.

As in any 2-category, two objects $\SSS,\tilde{\SSS}$ in $\mathbf{RigCat}$ (or in $\mathbf{SRigCat}$) are said to be equivalent when there exists (symmetric) rig category homomorphisms $\FF:\SSS\to\tilde{\SSS}$ and $\tilde{\FF}:\tilde{\SSS}\to\SSS$ and invertible rig transformations $\xi:\tilde{\FF}\circ\FF\To id_\SSS$ and $\tilde{\xi}:\FF\circ\tilde{\FF}\To id_{\tilde{\SSS}}$.

\subsection{Strictification theorem}
A generic (symmetric) rig category involves many natural isomorphisms and lots of required commutative diagrams. Hence it is useful to know that some of the natural isomorphisms can be assumed to be identities because the final structure is equivalent to a similar one but with some of these isomorphisms trivial. Theorems of this type are usually known as strictification theorems. For symmetric rig categories the theorem is due to May (\cite{May-1977}, Proposition~VI.3.5), and for generic rig categories it is a consequence of the more general strictification theorem for {\bf SMC}-categories due to Guillou \cite{Guillou-2010}. Their statements are as follows.

\subsubsection{\sc Theorem.}(\cite{May-1977},\cite{Guillou-2010})\label{teorema_estrictificacio}
Every rig category (resp. symmetric rig category) $\SSS$ is equivalent in $\mathbf{RigCat}$ (resp. in $\mathbf{SRigCat}$) to a semistrict rig category (resp. semistrict symmetric rig category).

\medskip
The choice of which distributor is made trivial in a semistrict version of a given (symmetric) rig category is logically arbitrary. The only relevant point here is that, in the absence of a strict commutativity of $+$, a common and usually unavoidable situation, it is unreasonable to demand that both distributors be identities.

\subsubsection{\sc Example} \label{Mat_k}
Let $k$ be any ground field $k$, and let $\Mm at_k$ be the category with objects the positive integers $n\geq 0$, with $0$ as zero object, and with the $m\times n$ matrices with entries in $k$ as morphisms $n\to m$ for each $m,n\geq 1$. When equipped with the sum and product given on objects in the usual way, and on morphisms by
\[
A+B=\left(\begin{array}{c|c} A&0 \\ \hline 0&B\end{array}\right),\qquad A\bullet B=\left(\begin{array}{ccc} A\,b_{11}&\cdots&A\,b_{1n} \\ \vdots&&\vdots \\ A\,b_{m1}&\cdots&A\,b_{mn}\end{array}\right)
\]
for any matrices $A,B$ with $B=(b_{ij})$, it becomes a left semistrict rig category $\MM at_k$. It provides a left semistrict (and skeletal) version of the rig category $\FF\VV ect_k$ of finite dimensional vector spaces over $k$.

\medskip
The details in the previous example can be omitted because they are not relevant in the sequel. By contrast, next section is devoted to a complete description of an explicit left semistrict version of the rig category of finite sets which is essential for what follows.

\section{Semistrict version of the symmetric 2-rig of finite sets}

The purpose of this paper is to show that the symmetric 2-rig $\widehat{\FF\SSS et}$ is biinitial in the 2-category of rig categories. However. instead of working with $\widehat{\FF\SSS et}$ we shall consider an equivalent, skeletal version of it we shall denote by $\widehat{\FF\SSS et}_{sk}$, which has the advantage of being semistrict. Since equivalent objects in a 2-category have equivalent categories of morphisms, it is indeed enough to prove that $\widehat{\FF\SSS et}_{sk}$ is biinitial in $\mathbf{RigCat}$. This considerably simplifies the diagrams, and makes computations much easier. 

This semistrict version of the $\widehat{\FF\SSS et}$ appears, for instance, as Example~VI.5.1 in \cite{May-1977}. However, the detailed description given here, including explicit descriptions of the distributors, seems to be new.

Let us first recall from \S~\ref{categories_distributives} that, after fixing particular binary products and coproducts, and final and initial objects every distributive category has a canonical symmetric rig category structure. In general, the resulting structural isomorphisms $a,c,l,r,a',c',l',r'$ are non trivial. In some cases, however, and for suitable choices of these binary products, coproducts and final, initial objects the associator and left and right unitors (but usually not the commutators) turn out to be trivial. This is so for the skeleton of $\Ff\Ss et$ having as objects the sets $[n]=\{1,\ldots,n\}$ for each $n\geq 1$, and $[0]=\emptyset$. We shall denote this skeleton by $\Ff\Ss et_{sk}$. 

\subsection{Additive monoidal structure}

Being a skeleton of $\Ff\Ss et$, the groupoid $\Ff\Ss et_{sk}$ has all binary products and coproducts, and $[1]$ and $[0]$ as (unique) final and initial objects, respectively. Moreover, chosing binary products and coproducts just amounts in this case to making appropriate choices of the respective projections and injections.

\subsubsection{\sc Lemma.}\label{estructura_additiva_semistricta}
For every objects $[m],[n]\in\Ff\Ss et_{sk}$ let us take as coproduct of the pair $([m],[n])$ the set $[m+n]$ with the injections $\iota^1_{[m],[n]}:[m]\to[m+n]$ and $\iota^2_{[m],[n]}:[n]\to[m+n]$ defined by 
\begin{align*}
\iota^1_{[m],[n]}(i)&=i, \quad i=1,\ldots m, 
\\ 
\iota^2_{[m],[n]}(j)&=m+j, \quad j=1,\ldots,n.
\end{align*}
Then for every maps $f:[m]\to[m']$, $g:[n]\to[n']$ their sum $f+g:[m+n]\to[m'+n']$ is given by
\begin{equation}\label{suma-morfismes}
(f+g)(k)=\left\{
\begin{array}{ll} 
f(k), & \mbox{if $k\in\{1,\ldots,m\}$,} 
\\
m'+g(k-m), & \mbox{if $k\in\{m+1,\ldots,m+n\}$,}
\end{array}
\right.
\end{equation}
and the resulting symmetric $+$-monoidal structure on $\Ff\Ss et_{sk}$ is semistrict with nontrivial commutators $c_{[m],[n]}:[m+n]\to[n+m]$ given by
\begin{equation}\label{commutador}
c_{[m],[n]}(k)=\left\{
\begin{array}{ll} 
n+k, & \mbox{if $k\in\{1,\ldots,m\}$,} 
\\
k-m, & \mbox{if $k\in\{m+1,\ldots,m+n\}$.}
\end{array}
\right.
\end{equation}
In particular, for each $n\geq 1$ the commutators $c_{[0],[n]}$ and $c_{[n],[0]}$ are both identities while $c_{[1],[n]}$ and $c_{[n],[1]}$ are the $(n+1)$-cycles $(1,n+1,n,n-1,\ldots,2)_{n+1}$ and $(1,2,3,\ldots,n+1)_{n+1}$, respectively.

\medskip
\noindent
{\sc Proof.} It basically follows from the fact that these injections are such that
\begin{eqnarray*}
&\iota^1_{[m+n],[p]}\,\iota^1_{[m],[n]}=\iota^1_{[m],[n+p]},
\\
&\iota^2_{[m],[n+p]}=\iota^2_{[m],[n]}+id_{[p]},
\\
&\iota^2_{[0],[n]}=id_{[n]}=\iota^1_{[n],[0]}
\end{eqnarray*}
for each $m,n,p\geq 0$. The details are left to the reader.
\qed

\subsubsection{\sc Remark.}
The injections $\iota^1_{[m],[n]},\iota^2_{[m],[n]}$ are nothing but the composite with the canonical injections of $[m],[n]$ into the disjoint union $[m]\sqcup[n]=([m]\times\{0\})\cup([n]\times\{1\})$ with the bijection $b^+_{m,n}:[m]\sqcup[n]\to[m+n]$ given by
\begin{equation}\label{b+}
b^+(k,\alpha)=\left\{ \begin{array}{ll}
k, & \mbox{if $\alpha=0$,} \\ m+k, & \mbox{if $\alpha=1$.}
\end{array}\right.
\end{equation}

\medskip
We shall denote by $\mathscr{F}\mathscr{S}et_{sk}^+$ (resp. $\widehat{\mathscr{F}\mathscr{S}et}_{sk}^+$) the semistrict (additive) symmetric monoidal category referred to in this lemma (resp. the underlying groupoid equipped with the inherited semistrict symmetric monoidal structure).

\subsection{Multiplicative monoidal structure}

The multiplicative monoidal structure on $\Ff\Ss et_{sk}$ is defined similarly. The object part of any product of $([m],[n])$ in $\Ff\Ss et_{sk}$ is necessarily the set $[mn]$, and the projections onto $[m]$ and $[n]$ can be obtained by chosing any natural bijection $b^{\,\bullet}_{m,n}:[mn]\to[m]\times[n]$, and composing it with the canonical projections. The point is that by suitably chosing the bijections $b^{\,\bullet}_{m,n}$ the final symmetric $\bullet$\,-monoidal structure is again semistrict.

\subsubsection{\sc Lemma.}\label{estructura_multiplicativa_semistricta}
For every objects $[m],[n]\in\Ff\Ss et_{sk}$ with $m,n\geq 1$ let us take as product of the pair $([m],[n])$ the set $[mn]$ with the projections $\pi^1_{[m],[n]}:[mn]\to[m]$ and $\pi^2_{[m],[n]}:[mn]\to[n]$ defined by 
\begin{align*}
\pi^1_{[m],[n]}(k)&=\left\{\begin{array}{ll} m, & \mbox{if $m\,|\,k$,} \\ r, & \mbox{otherwise,} \end{array}\right.
\\ 
\pi^2_{[m],[n]}(k)&=\left\{\begin{array}{ll} q, & \mbox{if $m\,|\,k$,} \\ q+1, & \mbox{otherwise,}\end{array}\right.
\end{align*}
where $q,r$ respectively denote the quotient and the remainder of the euclidean division of $k$ by $m$. Then for every maps $f:[m]\to[m']$, $g:[n]\to[n']$ with $m,m',n,n'\geq 1$, the product map $f\bullet g:[mn]\to[m'n']$ is given by
\begin{equation}\label{producte-morfismes}
(f\bullet g)(k)=\left\{
\begin{array}{ll} 
(g(q)-1)\,m'+f(m), & \mbox{if $m\,|\,k$,} 
\\
(g(q+1)-1)\,m'+f(r), & \mbox{otherwise,}
\end{array}
\right.
\end{equation}
and the resulting symmetric $\bullet$\,-monoidal structure on $\Ff\Ss et_{sk}$ is semistrict with nontrivial commutators $c'_{[m],[n]}:[mn]\to[nm]$, $m,n\geq 1$, given by
\begin{equation}\label{commutador'}
c'_{[m],[n]}(k)=\left\{
\begin{array}{ll} 
(m-1)n+q, & \mbox{if $m\,|\,k$,} 
\\
(r-1)n+q+1, & \mbox{otherwise,}
\end{array}
\right.
\end{equation}
where $q,r$ are as before. In particular, for each $n\geq 1$ the commutators $c'_{[1],[n]},c'_{[n],[1]}$ are both identities while $c'_{[2],[n]}$ and $c'_{[n],[2]}$ are given by
\begin{align*}
c'_{[2],[n]}(k)&=\left\{
\begin{array}{ll}
n+(k/2), & \mbox{if $k$ is even,} \\ (k+1)/2, & \mbox{if $k$ is odd,} 
\end{array}\right.
\\
c'_{[n],[2]}(k)&=\left\{
\begin{array}{ll}
2k-1, & \mbox{if $k\in\{1,\ldots,n\}$,} \\
2(k-n), & \mbox{if $k\in\{n+1,\ldots,2n\}$.}  
\end{array}\right.
\end{align*}

\medskip
\noindent
{\sc Proof.} For each $m,n\geq 1$ let $b^{\,\bullet}_{m,n}:[m]\times[n]\to[mn]$ be the bijection given by enumerating the elements of $[m]\times[n]$ by rows, i.e.
\[
b^{\,\bullet}_{m,n}(i,j)=(j-1)\,m+i.
\]
Its inverse $(b^{\,\bullet})^{-1}_{m,n}:[mn]\to[m]\times[n]$ is given by
\[
(b^{\,\bullet})^{-1}_{m,n}(k)=\left\{
\begin{array}{ll} (m,q), & \mbox{if $m\,|\,k$,} \\ (r,q+1), & \mbox{otherwise,} \end{array}\right.
\]
with $q,r$ as in the statement. Then the projections $\pi^1_{[m],[n]},\pi^2_{[m],[n]}$ make the diagram
\[
\xymatrix @R=.8pc @C=.8pc{
& [m]\times[n]\ar[d]^{b^{\,\bullet}_{m,n}}\ar[dl]_{pr_1}\ar[dr]^{pr_2} &
\\
[m] & [mn]\ar[l]^{\pi^1_{[m],[n]}}\ar[r]_{\pi^2_{[m],[n]}} & [n]
}
\]
commute, and (\ref{producte-morfismes}) is nothing but the composite
\[
f\bullet g=b^{\,\bullet}_{m',n'}\,(f\times g)\, (b^{\,\bullet})^{-1}_{m,n},
\]
with $f\times g$ the usual cartesian product map given by $(i,j)\mapsto(f(i),g(j))$. It follows that the diagram 
\[
\xymatrix @R=.7pc @C=.9pc{
[m]\ar[ddd]_f 
&& 
[mn]\ar[ll]_{\pi^1_{[m],[n]}}\ar[rr]^{\pi^2_{[m],[n]}}\ar[ddd]^{f\bullet\, g}
&&
[n]\ar[ddd]^g
\\
& [m]\times[n]\ar[ru]_{b^{\,\bullet}_{m,n}}\ar[lu]^{pr_1}\ar[d]^{f\times g} & & [m]\times[n]\ar[lu]^{b^{\,\bullet}_{m,n}}\ar[ru]_{pr_2}\ar[d]^{f\times g} &
\\
& [m']\times[n']\ar[ld]_{pr_1}\ar[rd]^{b^{\,\bullet}_{m',n'}} && [m']\times[n']\ar[rd]^{pr_2}\ar[ld]_{b^{\,\bullet}_{m',n'}}&
\\
[m'] && [m'n']\ar[ll]^{\pi^1_{[m'],[n']}}\ar[rr]_{\pi^2_{[m'],[n']}} && [n']   
}
\]
commutes, and hence (\ref{producte-morfismes}) makes the dual of the diagram in (D1) commute for the projections $\pi^1_{[m],[n]},\pi^2_{[m],[n]}$ in the statement. Moreover, these projections are such that
\begin{eqnarray*}
&\pi^1_{[m],[n]}\,\pi^1_{[mn],[p]}=\pi^1_{[m],[np]},
\\
&\pi^2_{[m],[np]}=\pi^2_{[m],[n]}\times id_{[p]},
\\
&\pi^2_{[1],[n]}=id_{[n]}=\pi^1_{[n],[1]}
\end{eqnarray*}
for each $m,n,p\geq 1$, which together with the uniqueness of the morphisms making the duals of the diagrams in (D3) and (D5) commute implies that all isomorphism $a'_{[m],[n],[p]},l'_{[m]},r'_{[m]}$ are identities. To prove these equalities, let be $l\in[mnp]$, and let
\[
l=q'(mn)+r'=q''m+r''
\]
be the euclidean divisions of $l$ by $mn$ and by $m$, respectively, and $r'=q'''m+r'''$ the euclidean division by $m$ of the remainder $r'$ of the first of these divisions. Clearly, we have
\begin{itemize}
\item[(i)] $r''=r'''$,
\item[(ii)] $q''=q'n+q'''$, and 
\item[(iii)] $m\,|\,l$ if and only if $m\,|\,r'$.
\end{itemize}
Then, on the one hand, an easy computation shows that
\[
(\pi^1_{[m],[n]}\circ\pi^1_{[mn],[p]})(l)=
\left\{\begin{array}{ll} m, & \mbox{if $mn\,|\,l$,} \\ m, & \mbox{if $mn\nmid l$ and $m\,|\,r'$,} \\ r''', & \mbox{if $mn\nmid l$ and $m\nmid r'$.} \end{array}\right.
\]
Since conditions $mn\,|\,l$ and $m\nmid r'$ can not hold simultaneously we have
\[
(\pi^1_{[m],[n]}\,\pi^1_{[mn],[p]})(l)=
\left\{\begin{array}{ll} m, & \mbox{if $m\,|\,r'$,} \\ r''', & \mbox{otherwise,} \end{array}\right.
\]
while by definition 
\[
\pi^1_{[m],[np]}(l)=
\left\{\begin{array}{ll} m, & \mbox{if $m\,|\,l$,} \\ r'', & \mbox{otherwise.} \end{array}\right.
\]
Hence it follows from (i) and (iii) that both maps are indeed the same. On the other hand, for each $l\in[mnp]$ we have
\begin{align*}
(\pi^2_{[m],[n]}\times id_{[p]})(l)&=
\left\{
\begin{array}{ll} (q'-1)n+\pi^2_{[m],[n]}(mn), & \mbox{if $mn\,|\,l$,} \\ q'n+\pi^2_{[m],[n]}(r'), & \mbox{if $mn\nmid l$,} \end{array}
\right.
\\
&=
\left\{
\begin{array}{ll} q'n, & \mbox{if $mn\,|\,l$,} \\ q'n+q''', & \mbox{if $mn\nmid l$ and $m\,|\,r'$,} \\ q'n+q'''+1, & \mbox{if $mn\nmid l$ and $m\nmid r'$,} \end{array}
\right.
\\
&=
\left\{
\begin{array}{ll} q'n+q''', & \mbox{if $m\,|\,r'$,} \\ q'n+q'''+1, & \mbox{otherwise} \end{array}
\right.
\end{align*}
(in the last equality we use that $mn\,|\,l$ implies that $r'=0$ and hence, $q'''=0$), while by definition 
\[
\pi^2_{[m],[np]}(l)=\left\{\begin{array}{ll} q'', & \mbox{if $m\,|\,l$,} \\ q''+1, & \mbox{otherwise.} \end{array}\right.
\]
Hence it follows from (ii) and (iii) that both maps are again the same. Finally, it readily follows from their definitions that $\pi^1_{[n],[1]}$ and $\pi^2_{[1],[n]}$ are identities. 

It remains to see that the map (\ref{commutador'}) makes the dual of the diagram in (D4) commute. Or, (\ref{commutador'}) is nothing but the composite
\[
c'_{[m],[n]}=b^{\,\bullet}_{n,m}\,\sigma_{[m],[n]}\, (b^{\,\bullet})^{-1}_{m,n},
\]
where $\tilde{c}\,'_{[m],[n]}:[m]\times[n]\to[n]\times[m]$ is the permutation map $(i,j)\mapsto(j,i)$. The commutativity of the dual of the diagram in (D4) follows then from the diagram
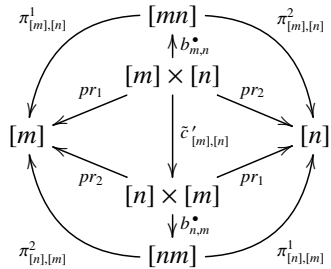
\[
\xymatrix @R=.7pc{
& [mn]\ar@/_1.5pc/[ldd]_{\pi^1_{[m],[n]}}\ar@/^1.5pc/[rdd]^{\pi^2_{[m],[n]}} & 
\\
& [m]\times[n]\ar[u]_{b^{\,\bullet}_{m,n}}\ar[ld]_{pr_1}\ar[rd]^{pr_2}\ar[dd]^{\tilde{c}\,'_{[m],[n]}} &
\\
[m] & & [n]
\\
& [n]\times[m]\ar[lu]^{pr_2}\ar[ur]_{pr_1}\ar[d]^{b^{\,\bullet}_{n,m}} & 
\\
& [nm]\ar@/^1.5pc/[luu]^{\pi^2_{[n],[m]}}\ar@/_1.5pc/[uur]_{\pi^1_{[n],[m]}} &
}
\]
all of whose inner triangles commute.
\qed

\subsubsection{\sc Remark.}
When we think of the elements in $[mn]$ as the points of the finite lattice $[m]\times [n]\subset\RR^2$, the map $f\bullet g$ simply corresponds to applying $f$ to the columns and $g$ to the rows. However, the explicit formula for $(f\bullet g)(k)$ depends on the way we decide to enumerate the points in the lattice and hence, on the chosen bijection $b^{\,\bullet}_{m,n}$. The same thing happens with the commutators, ultimately defined by the maps $(i,j)\mapsto(j,i)$. The above bijections $b^{\,\bullet}_{m,n}$ correspond to enumerating the points by rows, so that the formula (\ref{commutador'}) for the commutators corresponds to doing the following. Take two sets of $mn$ aligned points, one on the top of the other. Divide the top set into $n$ boxes each one with $m$ points, and the bottom set into $m$ boxes each one with $n$ points. Then $c'_{[m],[n]}$ maps the successive points in the top $j$ box, for each $j\in\{1,\ldots,n\}$, into the $j^{th}$ point in the successive $m$ bottom boxes.   

\subsection{Distributors and absorbing isomorphisms}
Next step is to describe the corresponding left and right distributors. Let $R_k:\NN\to[k]$ be the modified remainder function mapping each nonnegative integer $x$ to the remainder of the euclidean division of $x$ by $k$ if $k\nmid x$, and to $k$ if $k\,|\,x$. Then we have the following.

\subsubsection{\sc Lemma.}\label{distribuidors}
Let $\Ff\Ss et_{sk}$ be equipped with the above additive and multiplicative semistrict symmetric monoidal structures. Then for every objects $[m],[n],[p]$ the left distributor $\overline{d}_{[m],[n],[p]}$ is the identity, while the right distributor $\overline{d}\,'_{[m],[n],[p]}:[mp+np]\to[mp+np]$ is the identity when some of the objects $[m],[n],[p]$ is $[0]$, and otherwise is given by
\[
\overline{d}\,'_{[m],[n],[p]}(s)=\left\{\begin{array}{ll} s+\displaystyle{\frac{s-R_m(s)}{m}\,n}, & \mbox{if $s\in\{1,\ldots,mp\}$}
\\
s-mp+\displaystyle{\frac{s-mp+n-R_n(s-mp)}{n}\,m}, & \mbox{if $s\in\{mp+1,\ldots,mp+np\}$}\end{array}\right.
\]
In particular, $\overline{d}\,'_{[m],[n],[p]}$ is the identity when some of the integers $m,n,p$ is zero, and $\overline{d}\,'_{[n],[1],[1]},\overline{d}\,'_{[1],[n],[1]}:[n+1]\to[n+1]$ are both identities, and $\overline{d}\,'_{[1],[1],[n]}:[2n]\to[2n]$ is the permutation given by
\begin{equation}\label{d'_1,1,n}
\overline{d}'_{[1],[1],[n]}(s)=\left\{\begin{array}{ll} 
2s-1, & \mbox{if $s\in\{1,\ldots,n\}$,} 
\\
2(s-n), & \mbox{if $s\in\{n+1,\ldots,2n\}$}
\end{array}\right.
\end{equation}
for each $n\geq 1$.

\medskip
\noindent
{\sc Proof.} Using (\ref{producte-morfismes}) it is easy to check that
\begin{align*}
id_{[m]}\times\iota^1_{[n],[p]}&=\iota^1_{mn,mn+mp},  \\ id_{[m]}\times\iota^2_{[n],[p]}&=\iota^2_{mp,mn+mp}.
\end{align*}
Hence $\overline{d}_{[m],[n],[p]}$ is the identity. As to the right distributor $\overline{d}\,'_{[m],[n],[p]}$, by definition we have
\[
\overline{d}\,'_{[m],[n],[p]}(s)=
\left\{\begin{array}{ll} (\iota^1_{[m],[n]}\times id_{[p]})(s), & \mbox{if $s\in\{1,\ldots,mp\}$,} \\ (\iota^2_{[m],[n]}\times id_{[p]})(s-mp), & \mbox{if $s\in\{mp+1,\ldots,mp+np\}$} \end{array}\right.
\]
for each $s\in[mp+np]$. Hence $\overline{d}\,'_{[m],[n],[p]}$ is the identity when some of the objects $[m],[n],[p]$ is $[0]$. Otherwise, by (\ref{producte-morfismes}) we have
\[
(\iota^1_{[m],[n]}\times id_{[p]})(k)=
\left\{\begin{array}{ll} k+(q-1)n, & \mbox{if $k=qm$,} \\
k+qn, & \mbox{if $k=qm+r$, $r\neq 0$,} \end{array}\right.
\]
for each $k\in[mp]$, and
\[
(\iota^2_{[m],[n]}\times id_{[p]})(l)=
\left\{\begin{array}{ll} l+q'm, & \mbox{if $l=q'n$,} \\
l+(q'+1)m, & \mbox{if $l=q'n+r'$, $r'\neq 0$,} \end{array}\right.
\]
for each $l\in[np]$. Hence
\[
\overline{d}\,'_{[m],[n],[p]}(s)=
\left\{\begin{array}{ll} s+(q-1)n, & \mbox{if $s\in\{1,\ldots,mp\}$ and $m\,|\,s$,} \\ s+qn, & \mbox{if $s\in\{1,\ldots,mp\}$ and $m\nmid s$,} 
\\
s+(q'-p)m, & \mbox{if $s\in\{mp+1,\ldots,mp+np\}$ and $n\,|\,(s-mp)$,} \\ s+(q'-p+1)m, & \mbox{if $s\in\{mp+1,\ldots,mp+np\}$ and $n\nmid (s-mp)$,}\end{array}\right.
\]
where $q$ is the quotient of the euclidean division of $s$ by $m$ for $s\leq mp$, and $q'$ the quotient of the euclidean division of $s-mp$ by $n$ when $s>mp$. It is left to the reader checking that, in terms of the modified remainder functions, this is the permutation in the statement.
\qed

\subsubsection{\sc Example.}
The first few non trivial distributors $\overline{d}\,'_{[1],[1],[n]}:[2n]\to[2n]$ are 
\begin{align*}
\overline{d}\,'_{[1],[1],[2]}&=(2,3)_4,
\\
\overline{d}\,'_{[1],[1],[3]}&=(2,3,5,4)_6,
\\
\overline{d}\,'_{[1],[1],[4]}&=(2,3,5)_8\,(4,7,6)_8,
\\
\overline{d}\,'_{[1],[1],[5]}&=(2,3,5,9,8,6)_{10}\,(4,7)_{10}.
\end{align*}
Their decompositions into disjoint cycles do not seem to follow an easy pattern. Hence, unlike the commutators, all of order two, the order of the distributors $\overline{d}\,'_{[1],[1],[n]}$ will be a non-trivial function of $n$.

\subsubsection{\sc Remark.}
There seems to be no obvious description of the generic right distributor $\overline{d}'_{[m],[n],[p]}$ as a composite of the generators of $S_{(m+n)p}$. However, such a description exists when $m=n=1$, and is given as follows. Thus we know from (\ref{d'_1,1,n}) that $(d'_{[1],[1],[p]})^{-1}$ maps $1,\ldots,p$ to the first $p$ odd positive integers, and $p+1,\ldots,2p$ to the first $p$ even positive integers. For each $i=2,\ldots,p$ let $\sigma_{i,p}\in S_{2p}$ be the permutation given by
\[
\sigma_{i,p}=(i,i+1)_{2p}\,(i+2,i+3)_{2p}\,\cdots (i+2(p-i),i+2(p-i)+1)_{2p}.
\]
Then the reader may easily chek that $\overline{d}'_{[1],[1],[p]}=\sigma_{2,p}\,\sigma_{3,p}\cdots \sigma_{p-1,p}\,\sigma_{p,p}$.

\subsubsection{\sc Theorem.}\label{versio_semistricta}
The symmetric rig category structure on $\Ff\Ss et_{sk}$ canonically associated to the choices of products and coproducts of Lemmas~\ref{estructura_additiva_semistricta} and \ref{estructura_multiplicativa_semistricta} is left semistrict. The left symmetric rig category so obtained  $\FF\SSS et_{sk}$ is equivalent to the symmetric rig category $\FF\SSS et$.

\smallskip
\noindent
{\sc Proof.}
The first assertion is a consequence of Proposition~\ref{estructura_categoria_rig_simetrica}, the previous three lemmas, and the fact that the projections $\pi^1_{[0],[n]}$ and $\pi^2_{[n],[0]}$ are in this case identities. As to the equivalence between $\FF\SSS et_{sk}$ and $\FF\SSS et$, an equivalence is given by the inclusion functor $\Jj:\Ff\Ss et_{sk}\hookrightarrow\Ff\Ss et$ with the additive and multiplicative symmetric monoidal structures defined by the above bijections $b^+_{m,n}$ and $b^{\,\bullet}_{m,n}$. 
\qed

\subsubsection{\sc Corollary}
$\widehat{\FF\SSS et}_{sk}$ is a left semistrict symmetric 2-rig equivalent to the symmetric 2-rig $\widehat{\FF\SSS et}$. 

\medskip
\noindent
{\sc Proof.} It follows the the previous theorem and Example~\ref{core}.
\qed

\subsubsection{\sc Remark.}
If instead of the above bijections $b^{\,\bullet}_{m,n}$ one uses those that correspond to enumerating the points of $[m]\times[n]$ by columns, the obtained symmetric rig category structure on $\Ff\Ss et_{sk}$ turns out to be {\em right} semistrict.

\section{Main theorem}

Let $\SSS$ be an arbitrary rig category, not necessarily symmetric, with $\zero$ and $\un$ as zero and unit objects, respectively. Our purpose is to prove that the category 
of rig category homomorphisms from $\widehat{\FF\SSS et}$ to $\SSS$ is equivalent to the terminal category. Since this category, up to equivalence, remains invariant when $\SSS$ is replaced by any rig category equivalent to $\SSS$,
 Theorem~\ref{teorema_estrictificacio} above allows us to assume without loss of generality that $\SSS$ is left semistrict. We shall do that for the rest of the paper.

\subsection{The basic maps of an arbitrary homomorphism}

Since the additive associator $a$ of $\SSS$ is trivial, for every object $x\in\Ss$ and every $n\geq 0$ there is an object $nx$ uniquely defined by
\[
nx=\left\{\begin{array}{ll} 
\zero, & \mbox{if $n=0$,} 
\\
x+\stackrel{n)}{\cdots}+x, & \mbox{if $n>0$.}
\end{array}\right.
\]
When $x=\un$ we shall write $\underline{n}$ instead of $n\un\in\Ss$. In particular, $\underline{0}=\zero$ and $\underline{1}=\un$.

Recall from Definition~\ref{morfisme_categories_semianell} that a homomorphism of rig categories $\FF:\widehat{\FF\SSS et}_{sk}\to\SSS$ is given by a symmetric $+$-monoidal functor $\mathscr{F}:\widehat{\mathscr{FS}et}^+_{sk}\to \mathscr{S}^+$, given by a triple $\mathscr{F}=(F,\varphi^+,\varepsilon^+)$, together with an additional $\bullet$-monoidal structure $(\varphi^\bullet,\varepsilon^\bullet)$ on $F$ satisfying the appropriate coherence conditions.
Fundamental for what follows is the next family of invertible maps associated to any such homomorphism $\FF$.

\subsubsection{\sc Definition}\label{morfismes_basics} Let $\FF:\widehat{\FF\SSS et}_{sk}\to\SSS$ be a homomorphism of rig categories. The {\em basic maps} of $\FF$ are the (invertible) maps $\tau_n:F[n]\to\underline{n}$, $n\geq 0$, defined by the recurrence relation
\begin{equation}\label{recurrencia_tau_1}
\tau_{n+1}=(\varepsilon^\bullet+\tau_n)\,\varphi^+_{[1],[n]},\quad n\geq 0,
\end{equation}
with the initial condition $\tau_0=\varepsilon^+$. 

\subsubsection{\sc Remark}\label{remarca_recurrencia_2}
For each $\FF:\widehat{\FF\SSS et}_{sk}\to\SSS$ as before, let $\tau'_n:F[n]\to\underline{n}$, $n\geq 0$, be the morphisms defined by the recurrence relation
\begin{equation}\label{recurrencia_tau_2}
\tau'_{n+1}=(\tau'_n+\varepsilon^\bullet)\,\varphi^+_{[n],[1]}, \quad n\geq 0,
\end{equation}
with the initial condition $\tau'_0=\varepsilon^+$. Then $\tau'_n=\tau_n$ for each $n\geq 0$. We proceed by induction on $n\geq 0$. By definition, we have $\tau'_0=\tau_0$. Let us assume that $\tau'_{k}=\tau_k$ for each $0\leq k\leq n$ for some $n\geq 0$. Then
\begin{align*}
\tau_{n+1}&=(\varepsilon^{\,\bullet}+\tau_{n})\,\varphi_{[1],[n]}^+
\\
&=(\varepsilon^{\,\bullet}+\tau'_{n})\,\varphi_{[1],[n]}^+
\\
&=(\varepsilon^{\,\bullet}+\tau'_{n})\,(id_{F[1]}+\varphi_{[n-1],[1]}^+)^{-1}\,(\varphi_{[1],[n-1]}^++id_{F[1]})\,\varphi_{[n],[1]}^+
\\
&=[\varepsilon^{\,\bullet}+\tau'_{n}(\varphi_{[n-1],[1]}^+)^{-1}]\,(\varphi_{[1],[n-1]}^++id_{F[1]})\,\varphi_{[n],[1]}^+
\\
&=(\varepsilon^{\,\bullet}+\tau'_{n-1}+\varepsilon^{\,\bullet})\,(\varphi_{[1],[n-1]}^++id_{F[1]})\,\varphi_{[n],[1]}^+
\\
&=[(\varepsilon^{\,\bullet}+\tau'_{n-1})\,\varphi_{[1],[n-1]}^++\varepsilon^{\,\bullet}]\,\varphi_{[n],[1]}^+
\\
&=[(\varepsilon^{\,\bullet}+\tau_{n-1})\,\varphi_{[1],[n-1]}^++\varepsilon^{\,\bullet}]\,\varphi_{[n],[1]}^+
\\
&=(\tau_{n}+\varepsilon^{\,\bullet})\,\varphi_{[n],[1]}^+
\\
&=(\tau'_{n}+\varepsilon^{\,\bullet})\,\varphi_{[n],[1]}^+
\\
&=\tau'_{n+1}.
\end{align*}
The third equality follows from the coherence axiom on the monoidality isomorphisms $\varphi^+$, and the remaining ones readily follow from the induction hypothesis and the functoriality of $+$.

\medskip
The importance of the maps $\tau_n$ come from the fact, shown below, that they completely determine the homomorphism $\FF$ once the action on objects of the underlying functor is given (cf. Section~\ref{conclusions}). Moreover, unlike the monoidality isomorphisms $\varphi^+_{[m][n]}$ and $\varphi^\bullet_{[m],[n]}$ for each $m,n\geq 0$, which must satisfy naturality and coherence conditions, the maps $\tau_n$ can be chosen in a completely arbitrary way. Next result shows how the additive and multiplicative monoidality isomorphisms of $\FF$ can be obtained from the basic maps. Later on, we shall see that the action on morphisms of the underlying functor $F$ is also completely given by the basic maps.

\subsubsection{\sc Proposition} \label{monoidalitat_tau} Let $(\tau_n)_{n\geq 0}$ be the basic maps of an arbitrary homomorphism of rig categories $\FF=(F,\varphi^+,\varepsilon^+,\varphi^\bullet,\varepsilon^\bullet)$ from $\widehat{\FF\SSS et}_{sk}$ to a left semistrict rig category $\SSS$. Then the additive and multiplicative monoidality isomorphisms are respectively given by
\begin{align}
    \varphi^+_{[m],[n]}&=(\tau_m+\tau_n)^{-1}\,\tau_{m+n}, \label{varphi+_tau}
    \\
    \varphi^\bullet_{[m],[n]}&=(\tau_m\bullet\tau_n)^{-1}\,\tau_{mn}\label{varphibullet_tau}
\end{align}
for each $m,n\geq 0$.

\medskip
\noindent
{\sc Proof.}
Let us prove (\ref{varphi+_tau}) or equivalently, that $\tau_{m+n}=(\tau_{m}+\tau_{n})\,\varphi^+_{[m],[n]}$
for each $m,n\geq 0$. We proceed by induction on $n\geq 0$ for any given $m\geq 0$. If $n=0$ this means checking that
\[
\tau_{m}=(\tau_{m}+id_{\underline{0}})\,(id_{F[m]}+\varepsilon^+)\,\varphi^+_{[m],[0]}
\]
for each $m\geq 0$, and this is true because $\SSS$ is assumed to be semistrict. Let us now assume that (\ref{varphi+_tau}) holds for some $n\geq 0$ and every $m\geq 0$. Then for every $m\geq 0$ we have
\begin{align*}
\tau_{m+n+1}&=(\tau_{m+n}+\varepsilon^{\,\bullet})\,\varphi^+_{[m+n],[1]}
\\
&=[(\tau_{m}+\tau_{n})\,\varphi^+_{[m],[n]}+\varepsilon^{\,\bullet}]\,\varphi^+_{[m+n],[1]}
\\
&=(\tau_{m}+\tau_{n}+\varepsilon^{\,\bullet})\,(\varphi^+_{[m],[n]}+id_{F[1]})\,\varphi^+_{[m+n],[1]}
\\
&=(\tau_{m}+\tau_{n}+\varepsilon^{\,\bullet})\,(id_{F[1]}+\varphi^+_{[n],[1]})\,\varphi^+_{[m],[n+1]}
\\
&=[\tau_{m}+(\tau_{n}+\varepsilon^{\,\bullet})\,\varphi^+_{[n],[1]}]\,\varphi^+_{[m],[n+1]}
\\
&=(\tau_{m}+\tau_{n+1})\,\varphi^+_{[m],[n+1]}.
\end{align*}
The first and the last equalities follow from Remark~\ref{remarca_recurrencia_2}, the second one by the induction hypothesis, the third and fifth ones by the functoriality of $+$, and the fourth one by the coherence axioms required on $\varphi^+$. The proof of (\ref{varphibullet_tau}) is completely similar but with $\bullet$ instead of $+$.
$\hfill\square$

\medskip
Although we will not need it, next result gives an explicit expression of the basic maps of $\FF$ in terms of the data defining $\FF$. 

\subsubsection{\sc Proposition}
For each homomorphism of rig categories $\FF:\widehat{\FF\SSS et}_{sk}\to\SSS$ the corresponding basic maps $\tau_n$ are given by $\tau_0=\varepsilon^+$, $\tau_1=\varepsilon^\bullet$, and
\[
\tau_n=(\varepsilon^\bullet+\stackrel{n)}{\cdots}+\varepsilon^\bullet)\,(id_{(n-2)F[1]}+\varphi^+_{[1],[1]})\,(id_{(n-3)F[1]}+\varphi^+_{[1],[2]})\,\cdots\,(id_{F[1]}+\varphi^+_{[1],[n-2]})\,\varphi^+_{[1],[n-1]}
\]
for each $n\geq 2$.

\medskip
\noindent
{\sc Proof.} The case $n=0$ is by definition of the basic maps. Moreover, when $n=0$ the recurrence relation (\ref{recurrencia_tau_1}) gives that
\[
\tau_1=(\varepsilon^\bullet+\varepsilon^+)\,\varphi^+_{[1],[0]}=(\varepsilon^\bullet+id_{\underline{0}})\,(id_{F[1]}+\varepsilon^+)\,\varphi^+_{[1],[0]}=\varepsilon^\bullet,
\]
where we have used that $\SSS$ is left semistrict (in particular, that the additive right unitor $\rho_{F[1]}$ is trivial) as well as the left semistrictness of $\widehat{\FF\SSS et}_{sk}$. The above expression for the remaining maps $\tau_n$ follows now by an easy induction on $n\geq 2$ using the functoriality of $+$. The details are left to the reader. 
$\hfill\square$

\subsection{The canonical homomorphism from finite sets to an arbitrary left semistrict rig category}
In order to see that, up to isomorphism, there is only one homomorphism of rig categories $\FF:\widehat{\FF\SSS et}_{sk}\to\SSS$ let us start by describing a particularly simple such homomorphism we shall call the canonical one.

\subsubsection{\sc Proposition}\label{morfisme_canonic}
Let $\SSS$ be an arbitrary left semistrict rig category, not necessarily symmetric. There exists a functor $F_{can}:\widehat{\Ff\Ss et}_{sk}\to\Ss$ acting on objects by $[n]\mapsto\underline{n}$, and on the generators of $S_n$ for each $n\geq 2$ by
\begin{equation}\label{accio_sobre_morfismes}
F_{can}((i,i+1)_n)=id_{\underline{i-1}}+c_{\underline{1},\underline{1}}+id_{\underline{n-i-1}}, \quad i=1,\ldots,n-1.
\end{equation}
Moreover, $F_{can}$ is the underlying functor of a {\em strict morphism} of rig categories $\FF_{can}:\widehat{\FF\SSS et}_{sk}\to\SSS$. 

\medskip
\noindent
{\sc Proof.} For short, let us denote by $\gamma_i$ the morphism in the right hand side of (\ref{accio_sobre_morfismes}). These assignments define a functor if and only if $\gamma_1,\ldots,\gamma_{n-1}$ satisfy for each $n\geq 2$ the same relations as the generators $(1,2)_n,\ldots,(n-1,n)_n$ of $S_n$, i.e. if and only if
\begin{itemize}
\item[(R1)] $\gamma_i^2=id_{\underline{n}}$ for each $i=1,\ldots,n-1$;
\item[(R2)] $(\gamma_i\,\gamma_{i+1})^3=id_{\overline{n}}$ for each $i=1,\ldots,n-2$;
\item[(R3)] $(\gamma_i\,\gamma_j)^2=id_{\overline{n}}$ for each $i,j=1,\ldots,n-1$ such that $|i-j|>1$.
\end{itemize} 
Relation (R1) follows from the functoriality of $+$, and the fact that $c$ is not just a braiding, but a symmetry, so that $c^2_{\underline{1},\underline{1}}=id_{\underline{2}}$. To see (R2), notice that by the functoriality of $+$ we have
\begin{align*}
\gamma_i\,\gamma_{i+1}&=(id_{\underline{i-1}}+c_{\underline{1},\underline{1}}+id_{\underline{n-i-1}})\,(id_{\underline{i}}+c_{\underline{1},\underline{1}}+id_{\underline{n-i-2}})
\\
&=(id_{\underline{i-1}}+c_{\underline{1},\underline{1}}+id_{\underline{1}}+id_{\underline{n-i-2}})\,(id_{\underline{i-1}}+id_{\underline{1}}+c_{\underline{1},\underline{1}}+id_{\underline{n-i-2}})
\\
&=id_{\underline{i-1}}+(c_{\underline{1},\underline{1}}+id_{\underline{1}})\,(id_{\underline{1}}+c_{\underline{1},\underline{1}})+id_{\underline{n-i-2}}).
\end{align*}
Hence (R2) holds if and only if
\[
[(c_{\underline{1},\underline{1}}+id_{\underline{1}})\,(id_{\underline{1}}+c_{\underline{1},\underline{1}})]^3=id_{\underline{3}}
\]
for each $i=1,\ldots,n-2$. Now, in every braided monoidal category $\mathscr{C}$ with trivial associator, and braiding $c$ the diagram
\[
\xymatrix{
x+y+z\ar[r]^{id_x+c_{y,z}}\ar[d]_{c_{x,y}+id_z} & x+z+y\ar[r]^{c_{x,z}+id_y} & z+x+y\ar[d]^{id_z+c_{x,y}} 
\\
y+x+z\ar[r]_{id_y+c_{x,z}} & y+z+x\ar[r]_{c_{y,z}+id_x} & z+y+x }
\]
commutes for every objects $x,y,z\in\Cc$ (cf. \cite{Joyal-Street-1993}, Proposition~2.1). In particular, when $\mathscr{C}=\mathscr{S}^+$, and $x=y=z=\underline{1}$ we obtain that
\[
(c_{\underline{1},\underline{1}}+id_{\underline{1}})\,(id_{\underline{1}}+c_{\underline{1},\underline{1}})=(id_{\underline{1}}+c_{\underline{1},\underline{1}})\,(c_{\underline{1},\underline{1}}+id_{\underline{1}})\,(id_{\underline{1}}+c_{\underline{1},\underline{1}})\,(c_{\underline{1},\underline{1}}+id_{\underline{1}}).
\]
The desired equality follows then by taking the composite of this equality on the left with $c_{\underline{1},\underline{1}}+id_{\underline{1}}$, and on the right with $id_{\underline{1}}+c_{\underline{1},\underline{1}}$. Finally, if $j\geq i+2$ we have
\begin{align*}
\gamma_i\,\gamma_j&=(id_{\underline{i-1}}+c_{\underline{1},\underline{1}}+id_{\underline{j-i}}+id_{\underline{n-j-1}})\,(id_{\underline{i}}+id_{\underline{j-i}}+c_{\underline{1},\underline{1}}+id_{\underline{n-j-1}})
\\
&=id_{\underline{i-1}}+(c_{\underline{1},\underline{1}}+id_{\underline{j-i}})\,(id_{\underline{j-i}}+c_{\underline{1},\underline{1}})+id _{\underline{n-j-1}}).
\end{align*}
Hence (R3) holds in this case if and only if
\[
[(c_{\underline{1},\underline{1}}+id_{\underline{j-i}})\,(id_{\underline{j-i}}+c_{\underline{1},\underline{1}})]^2=id_{j-i+2}.
\]
Now, when $j-i=2$ this is equal to the square of $c_{\underline{1},\underline{1}}+c_{\underline{1},\underline{1}}$ and hence, $id_{\underline{4}}$ because $c$ is a symmetry, and when $j-i>2$ it is the square of
\[
(c_{\underline{1},\underline{1}}+id_{\underline{j-i-2}}+id_{\underline{2}})\,(id_{\underline{2}}+id_{\underline{j-i-2}}+c_{\underline{1},\underline{1}})=c_{\underline{1},\underline{1}}+id_{\underline{j-i-2}}+c_{\underline{1},\underline{1}},
\]
and hence, $id_{\underline{j-i+2}}$. The case $i\geq j+2$ is proved similarly.

Clearly, the functor $F_{can}$ so defined preserves strictly the zero and unit objects. It further preserves the action of $+$ and $\bullet$ on objects strictly. Thus, being $\SSS$ semistrict, the additive associator is trivial and hence, for every $m,n\geq 0$ we have
\[
F_{can}([m]+[n])=F_{can}([m+n])=\underline{m+n}=\underline{m}+\underline{n}=F_{can}([m])+F_{can}([n]).
\]
Moreover, the triviality of the left distributors and of the multiplicative right unitors also gives that
\begin{align*}
\underline{m}\bullet\underline{n}&=\underline{m}\bullet(1+\stackrel{n)}{\cdots}+1)
\\
&=\underline{m}\bullet 1+\stackrel{n)}{\cdots}+\underline{m}\bullet 1
\\
&=\underline{m}+\stackrel{n)}{\cdots}+\underline{m}
\\
&=\underline{mn},
\end{align*}
and hence
\[
F_{can}([m]\bullet [n])=F_{can}([mn])=\underline{mn}=\underline{m}\bullet\underline{n}=F_{can}([m])\bullet F_{can}([n]).
\]
Therefore in order to see that $F_{can}$ is the underlying functor of a strict morphism of rig categories from $\widehat{\FF\SSS} et_{sk}$ to $\SSS$ it remains to be checked that $F_{can}$ also preserves strictly:
\begin{itemize}
\item[(SH1)] the action of $+$ and $\bullet$ on morphisms, 
\item[(SH2)] the additive commutators, and
\item[(SH3)] the right distributors.
\end{itemize}
Notice that the remaining structural isomorphisms are automatically preserved because $\widehat{\FF\SSS et}_{sk}$ is left semistrict, and $\SSS$ is assumed to be so.

\smallskip
\noindent
\underline{Proof of (SH1)}. For every $m\geq 1$ and $n\geq 2$ we have
\begin{align*}
F_{can}(id_{[m]}+(j,j+1)_n)&=F_{can}((m+j,m+j+1)_{m+n})
\\
&=id_{\underline{m+j-1}}+c_{\underline{1},\underline{1}}+id_{\underline{n-j-1}}
\\
&=id_{\underline{m}}+id_{\underline{j-1}}+c_{\underline{1},\underline{1}}+id_{\underline{n-j-1}}
\\
&=F_{can}(id_{[m]})+F_{can}((j,j+1)_n)
\end{align*}
for each $j\in\{1,\ldots,n-1\}$. Since both $+$ and $F_{can}$ are functorial, it follows that for every $\sigma\in S_n$ with $\sigma=(j_1,i_1+1)_n,\cdots(j_l,j_l+1)_n$ we have
\begin{align*}
F_{can}(id_{[m]}+\sigma)&=F_{can}((id_{[m]}+(j_1,i_1+1)_n)\cdots(id_{[m]}+(j_l,j_l+1)_n))
\\
&=F_{can}(id_{[m]}+(j_1,i_1+1)_n)\cdots F_{can}(id_{[m]}+(j_l,j_l+1)_n))
\\
&=(id_{\underline{m}}+F_{can}((j_1,j_1+1)_n))\cdots(id_{\underline{m}}+F_{can}((j_l,j_l+1)_n))
\\
&=id_{\underline{m}}+F_{can}((j_1,j_1+1)_n)\cdots F_{can}((j_l,j_l+1)_n)
\\
&=id_{\underline{m}}+F_{can}((j_1,j_1+1)_n\cdots (j_l,j_l+1)_n)
\\
&=id_{\underline{m}}+F_{can}(\sigma)
\end{align*}
Moreover, for every $m\geq 2$ and $n\geq 1$ we have
\begin{align*}
F_{can}((i,i+1)_m+id_{[n]})&=F_{can}((i,i+1)_{m+n})
\\
&=id_{\underline{i-1}}+c_{\underline{1},\underline{1}}+id_{\underline{m+n-i-1}}
\\
&=id_{\underline{i-1}}+c_{\underline{1},\underline{1}}+id_{\underline{m-i-1}}+id_{\underline{n}}
\\
&=F_{can}((i,i+1)_m)+F_{can}(id_{[n]})
\end{align*}
for each $i\in\{1,\ldots,m-1\}$ and the same argument as before shows that
\[
F_{can}(\rho+id_{[n]})=F_{can}(\rho)+id_{\underline{n}}
\]
for every $\rho\in S_m$. Using again the functoriality of $+$ and $F_{can}$, it follows that
\begin{align*}
F_{can}(\rho+\sigma)&=F_{can}((\rho+id_{[n]})\,(id_{[m]}+\sigma))
\\
&=F_{can}(\rho+id_{[n]})\,F_{can}(id_{[m]}+\sigma)
\\
&=(F_{can}(\rho)+id_{\underline{n}})\,(id_{\underline{m}}+F_{can}(\sigma))
\\
&=F_{can}(\rho)+F_{can}(\sigma)
\end{align*}
 for every permutations $\rho\in S_m$ and $\sigma\in S_n$. Hence $F_{can}$ preserves strictly the action of $+$ on morphisms.

To prove that it also preserves strictly the action of $\bullet$\,, we proceed in the same way. Thus it is enough to prove the special cases
\begin{itemize}
\item[(i)] $F_{can}((i,i+1)_m\bullet id_{[n]})=F_{can}((i,i+1)_m)\bullet id_{\underline{n}}$ for each $i=1,\ldots,m-1$ and $m\geq 2$, and
\item[(ii)] $F_{can}(id_{[m]}\bullet(j,j+1)_n)=id_{\underline{m}}\bullet F_{can}((j,j+1)_n)$ for each $j=1,\ldots,n-1$ and $n\geq 2$.
\end{itemize}
The generic case follows then by the functoriality of $F_{can}$ and $\bullet$\,. Let us prove (i). On the one hand, by definition of $\bullet$ we have
\[
(i,i+1)_m\bullet id_{[n]}=(i,i+1)_{mn}\,(m+i,m+i+1)_{mn}\,(2m+i,2m+i+1)_{mn}\cdots ((n-1)m+i,(n-1)m+i+1)_{mn}.
\]
Hence
\begin{align*}
F_{can}&((i,i+1)_m\bullet id_{[n]})
\\
&=F_{can}((i,i+1)_{mn})\, F_{can}((m+i,m+i+1)_{mn})\, F_{can}((2m+i,2m+i+1)_{mn})
\\
&\hspace{8truecm}\cdots F_{can}(((n-1)m+i,(n-1)m+i+1)_{mn})
\\
&=(id_{\underline{i-1}}+c_{\underline{1},\underline{1}}+id_{\underline{nm-i-1}})\,(id_{\underline{m+i-1}}+c_{\underline{1},\underline{1}}+id_{\underline{(n-1)m-i-1}})\,(id_{\underline{2m+i-1}}+c_{\underline{1},\underline{1}}+id_{\underline{(n-2)m-i-1}})
\\
&\hspace{9truecm}\cdots 
(id_{\underline{(n-1)m+i-1}}+c_{\underline{1},\underline{1}}+id_{\underline{m-i-1}})
\\
&=id_{\underline{i-1}}+(c_{\underline{1},\underline{1}}+id_{\underline{(n-1)m}})\,(id_{\underline{m}}+c_{\underline{1},\underline{1}}+id_{\underline{(n-2)m}})\,(id_{\underline{2m}}+c_{\underline{1},\underline{1}}+id_{\underline{(n-3)m}})
\cdots (id_{\underline{(n-1)m}}+c_{\underline{1},\underline{1}})+id_{\underline{m-i-1}},
\end{align*}
where in the last equality we use the functoriality of $+$ in order to separate the summand $id_{\underline{i-1}}$ on the left, and the summand $id_{\underline{m-i-1}}$ on the right which are common to all factors. Now, since $m\geq 2$ the composition of permutations in the middle is equal to
\begin{align*}
(c_{\underline{1},\underline{1}}+id_{\underline{(n-1)m}})\,(id_{\underline{m}}+c_{\underline{1},\underline{1}}+id_{\underline{(n-2)m}})\,(id_{\underline{2m}}+&c_{\underline{1},\underline{1}}+id_{\underline{(n-3)m}})\cdots (id_{\underline{(n-1)m}}+c_{\underline{1},\underline{1}})
\\
&=c_{\underline{1},\underline{1}}+id_{\underline{m-2}}+c_{\underline{1},\underline{1}}+id_{\underline{m-2}}+\cdots +id_{\underline{m-2}}+c_{\underline{1},\underline{1}},
\end{align*}
where we have $n$ summands equal to $c_{\underline{1},\underline{1}}$, and $n-1$ summands equal to $id_{\underline{m-2}}$. Hence
\[
F_{can}((i,i+1)_m\bullet id_{[n]})=id_{\underline{i-1}}+c_{\underline{1},\underline{1}}+id_{\underline{m-2}}+c_{\underline{1},\underline{1}}+id_{\underline{m-2}}+\cdots +id_{\underline{m-2}}+c_{\underline{1},\underline{1}}+id_{\underline{m-i-1}}.
\]
On the other hand, since $\SSS$ is assumed to be left semistrict, the left distributors and the multiplicative right unitors are trivial. Therefore
\begin{align*}
F_{can}((i,i+1)_m)\bullet id_{\underline{n}}&=(id_{\underline{i-1}}+c_{\underline{1},\underline{1}}+id_{\underline{m-i-1}})\bullet id_{\underline{n}}
\\
&=(id_{\underline{i-1}}+c_{\underline{1},\underline{1}}+id_{\underline{m-i-1}})\bullet (id_{\underline{1}}+\stackrel{n)}{\cdots}+id_{\underline{1}})
\\
&=(id_{\underline{i-1}}+c_{\underline{1},\underline{1}}+id_{\underline{m-i-1}})+\stackrel{n)}{\cdots}+(id_{\underline{i-1}}+c_{\underline{1},\underline{1}}+id_{\underline{m-i-1}})
\\
&=id_{\underline{i-1}}+c_{\underline{1},\underline{1}}+id_{\underline{m-2}}+c_{\underline{1},\underline{1}}+id_{\underline{m-2}}+\cdots +id_{\underline{m-2}}+c_{\underline{1},\underline{1}}+id_{\underline{m-i-1}},
\end{align*}
where the number of summands equal to $c_{\underline{1},\underline{1}}$ and to $id_{\underline{m-2}}$ are as before. This proves (i). Let us now prove (ii). By definition of $\bullet$\,, we have
\begin{equation}\label{descomposicio}
id_{[m]}\bullet(j,j+1)_n=((j-1)m+1,jm+1)_{mn}\,((j-1)m+2,jm+2)_{mn}\cdots((j-1)m+m,jm+m)_{mn}.
\end{equation}
Unlike before, the images by $F_{can}$ of the transpositions $((j-1)m+i,jm+i)_{mn}$ for each $i=1,\ldots,m$ are not immediate because they are not part of the generators of $S_{mn}$, whose images are given by Eq.~(\ref{accio_sobre_morfismes}). In fact, to compute their images we need to decompose $((j-1)m+i,jm+i)_{mn}$ as a product of the generators of $S_{mn}$. As the reader may easily check, such a decomposition is given by
\begin{align*}
((j-1)m+i,jm+i)_{mn}&=(jm+i-1,jm+i)_{mn}\,(jm+i-2,jm+i-1)_{mn}
\\
&\hspace{1truecm}\cdots\,(jm+i-(m-1),jm+i-(m-2))_{mn}\,(jm+i-m,jm+i-(m-1))_{mn}.
\end{align*}
Hence
\begin{align*}
F_{can}(((j-1)&m+i,jm+i)_{mn})
\\
&=F_{can}((jm+i-1,jm+i)_{mn})\,F_{can}((jm+i-2,jm+i-1)_{mn})
\\
&\hspace{0.8truecm}\cdots\,F_{can}((jm+i-(m-1),jm+i-(m-2))_{mn})\,F_{can}((jm+i-m,jm+i-(m-1))_{mn})
\\
&=(id_{\underline{jm+i-2}}+c_{\underline{1},\underline{1}}+id_{\underline{m(n-j)-i}})\,(id_{\underline{jm+i-3}}+c_{\underline{1},\underline{1}}+id_{\underline{m(n-j)-i+1}})
\\
&\hspace{0.8truecm}\cdots (id_{\underline{jm+i-m}}+c_{\underline{1},\underline{1}}+id_{\underline{m(n-j)-i+(m-2)}})\,(id_{\underline{jm+i-(m+1)}}+c_{\underline{1},\underline{1}}+id_{\underline{m(n-j)-i+(m-1)}})
\\
&=id_{\underline{jm+i-(m+1)}}+(id_{\underline{m-1}}+c_{\underline{1},\underline{1}})\,(id_{\underline{m-2}}+c_{\underline{1},\underline{1}}+id_{\underline{1}})
\\
&\hspace{2.2truecm}\cdots (id_{\underline{1}}+c_{\underline{1},\underline{1}}+id_{\underline{m-2}})\,(c_{\underline{1},\underline{1}}+id_{\underline{m-1}})+id_{\underline{m(n-j)-i}},
\end{align*}
where in the last equality we again make use of the functoriality of $+$ in order to separate what is common to all factors, namely, the summand $id_{\underline{jm+i-(m+1)}}$ on the left, and the summand $id_{\underline{m(n-j)-i}}$ on the right. An easy induction on $m\geq 1$ shows that the composition of permutations in the middle is just $c_{\underline{1},\underline{m}}$. When $m=1$ the composition indeed reduces to $c_{\underline{1},\underline{1}}$. Let us now assume that
\[
(id_{\underline{m-1}}+c_{\underline{1},\underline{1}})\,(id_{\underline{m-2}}+c_{\underline{1},\underline{1}}+id_{\underline{1}})\cdots (id_{\underline{1}}+c_{\underline{1},\underline{1}}+id_{\underline{m-2}})\,(c_{\underline{1},\underline{1}}+id_{\underline{m-1}})=c_{\underline{1},\underline{m}}
\]
for some $m\geq 1$. Then
\begin{align*}
(id_{\underline{m}}+c_{\underline{1},\underline{1}})\,&(id_{\underline{m-1}}+c_{\underline{1},\underline{1}}+id_{\underline{1}})\,(id_{\underline{m-2}}+c_{\underline{1},\underline{1}}+id_{\underline{2}})\cdots (id_{\underline{1}}+c_{\underline{1},\underline{1}}+id_{\underline{m-1}})\,(c_{\underline{1},\underline{1}}+id_{\underline{m}})
\\
&=(id_{\underline{m}}+c_{\underline{1},\underline{1}})\,[(id_{\underline{m-1}}+c_{\underline{1},\underline{1}})\,(id_{\underline{m-2}}+c_{\underline{1},\underline{1}}+id_{\underline{1}})\cdots (id_{\underline{1}}+c_{\underline{1},\underline{1}}+id_{\underline{m-2}})\,(c_{\underline{1},\underline{1}}+id_{\underline{m-1}})+id_{\underline{1}}]
\\
&=(id_{\underline{m}}+c_{\underline{1},\underline{1}})\,(c_{\underline{1},\underline{m}}+id_{\underline{1}})
\\
&=c_{\underline{1},\underline{m+1}},
\end{align*}
where in the first equality we have separated the right summand $id_{\underline{1}}$ common to all but the first factor, and in the second we apply the induction hypothesis. As to the third equality, it follows from the hexagon coherence axiom on $c$, according to which in any braided (in particular, symmetric) monoidal category $\mathscr{C}$ with trivial associator and braiding $c$, the diagram
\begin{equation}\label{axioma_hexagon}
\xymatrix{
x\otimes y\otimes z\ar[r]^{c_{x,y}\otimes id_z}\ar[rd]_{c_{x,y\otimes z}} & y\otimes x\otimes z\ar[d]^{id_y\otimes c_{x,z}}
\\
& y\otimes z\otimes x}
\end{equation}
commutes for every objects $x,y,z\in\Cc$. Thus we conclude that
\[
F_{can}(((j-1)m+i,jm+i)_{mn})=id_{\underline{jm+i-(m+1)}}+c_{\underline{1},\underline{m}}+id_{\underline{m(n-j)-i}}
\]
for each $i=1,\ldots,m$. Coming back to Eq.~(\ref{descomposicio}), it follows that
\begin{align*}
F_{can}(id_{[m]}\bullet&(j,j+1)_n)
\\
&=F_{can}(((j-1)m+1,jm+1)_{mn})\,F_{can}(((j-1)m+2,jm+2)_{mn})
\\ &\hspace{1.5truecm}\cdots\,F_{can}(((j-1)m+(m-1),(j-1)m+(m-1))_{mn})\,F_{can}(((j-1)m+m,jm+m)_{mn})
\\
&=(id_{\underline{jm-m}}+c_{\underline{1},\underline{m}}+id_{\underline{m(n-j)-1}})\,(id_{\underline{jm-m+1}}+c_{\underline{1},\underline{m}}+id_{\underline{m(n-j)-2}})
\\ &\hspace{2truecm}\cdots\,(id_{\underline{jm-2}}+c_{\underline{1},\underline{m}}+id_{\underline{m(n-j)-(m-1)}})\,(id_{\underline{jm-1}}+c_{\underline{1},\underline{m}}+id_{\underline{m(n-j)-m}})
\\
&=id_{\underline{jm-m}}+(c_{\underline{1},\underline{m}}+id_{\underline{m-1}})\,(id_{\underline{1}}+c_{\underline{1},\underline{m}}+id_{\underline{m-2}})\,
(id_{\underline{m-2}}+c_{\underline{1},\underline{m}}+id_{\underline{1}})\,(id_{\underline{m-1}}+c_{\underline{1},\underline{m}})+id_{\underline{m(n-j)-m}}.
\end{align*}
It is now easy to show by induction on $m\geq 1$, using again the above hexagon axiom but with all arrows reversed, that the composition of permutations in the middle is nothing but $c_{\underline{m},\underline{m}}$. In summary, we have
\[
F_{can}(id_{[m]}\bullet(j,j+1)_n)=id_{\underline{jm-m}}+c_{\underline{m},\underline{m}}+id_{\underline{m(n-j)-m}}.
\]
As to the right hand side of (ii), since $\SSS$ is assumed to be left semistrict, we have
\begin{align*}
id_{\underline{m}}\bullet F_{can}((j,j+1)_n)&=id_{\underline{m}}\bullet (id_{\underline{j-1}}+c_{\underline{1},\underline{1}}+id_{\underline{n-j-1}})
\\
&=id_{\underline{jm-m}}+id_{\underline{m}}\bullet c_{\underline{1},\underline{1}}+id_{\underline{m(n-j)-m}}
\end{align*}
Thus in order to prove (ii) it is now enough to see that
\[
id_{\underline{m}}\bullet c_{\underline{1},\underline{1}}=c_{\underline{m},\underline{m}},
\]
and this equality follows from the fact that the left distributor of $\SSS$ is trivial, together with the axiom
\[
\xymatrix{
x(y+z)\ar[r]^{id_x\bullet c_{y,z}}\ar[d]_{d_{x,y,z}} &x(z+y)\ar[d]^{d_{x,z,y}} 
\\
xy+xz\ar[r]_{c_{xy,xz}} & xz+xy}
\]
which holds in every rig category. Therefore $F_{can}$ also preserves strictly the action of $\bullet$ on morphisms.

\bigskip
\noindent
\underline{Proof of (SH2)}. We need to prove that $F_{can}(c_{[m],[n]})=c_{\underline{m},\underline{n}}$ for each $m,n\geq 0$. The cases $m=0$ or $n=0$ follow from the fact that in every semistrict rig category $\SSS$ the commutators $c_{0,x}$ and hence, also $c_{x,0}=c^{-1}_{0,x}$ are identities for every object $x\in\Ss$. For instance, if $\SSS$ is left semistrict, the left distributors and the absorbing isomorphisms are all trivial and hence, the previous commutative diagram with $y=0$ gives that
\[
id_x\bullet c_{0,z}=c_{0,xz}
\]
for every $x,z\in\Ss$. In particular, if $z=1$ we obtain that
\[
id_x\bullet c_{0,1}=c_{0,x}
\]
because the right multiplicative unitors are also identities. Moreover, $c_{0,1}=id_1$ because of the coherence axiom $r_1\,c_{0,1}=l_1$ and the triviality of the additive unitors. When $\SSS$ is right semistrict, a similar argument works taking as starting point the axiom
\[
\xymatrix{
(y+z)x\ar[r]^{c_{y,z}\bullet\, id_x}\ar[d]_{d'_{y,z,x}} &(z+y)x\ar[d]^{d'_{z,y,x}} 
\\
yx+zx\ar[r]_{c_{yx,zx}} & zx+yx}
\]
which also holds in every rig category. Since there is no obvious description of the commutators $c_{[m],[n]}$ in terms of the generators of $S_{m+n}$ except when $m=1$ or $n=1$, the cases $m,n\geq 1$ are shown by induction on $n\geq 1$ for any given $m\geq 1$. For $n=1$ we know that
\[
c_{[m],[1]}=(1,2,\ldots,m+1)_{m+1}=(1,2)_{m+1}\,(2,3)_{m+1}\cdots (m,m+1)_{m+1}.
\]
Hence
\begin{align*}
F_{can}(c_{[m],[1]})&=F_{can}((1,2)_{m+1})\,F_{can}((2,3)_{m+1})\cdots F_{can}((m,m+1)_{m+1})
\\
&=(c_{\underline{1},\underline{1}}+id_{\underline{m-1}})\,(id_{\underline{1}}+c_{\underline{1},\underline{1}}+id_{\underline{m-2}})\cdots (id_{\underline{m-1}}+c_{\underline{1},\underline{1}})
\\
&=c_{\underline{m},\underline{1}}.
\end{align*}
The last equality follows by an easy induction on $m\geq 1$ using again the hexagon axiom on the commutator. Let now be $n\geq 2$, and let us assume that $F_{can}(c_{[m],[n-1]})=c_{\underline{m},\underline{n-1}}$. Because of the hexagon axiom we have
\[
c_{[m],[n]}=(id_{[n-1]}+c_{[m],[1]})\,(c_{[m],[n-1]}+id_{[1]}).
\]
Therefore
\begin{align*}
F_{can}(c_{[m],[n]})&=F_{can}(id_{[n-1]}+c_{[m],[1]})\,F_{can}(c_{[m],[n-1]}+id_{[1]})
\\
&=(id_{\underline{n-1}}+F_{can}(c_{[m],[1]}))\,(F_{can}(c_{[m],[n-1]})+id_{\underline{1}})
\\
&=(id_{\underline{n-1}}+c_{\underline{m},\underline{1}})\,(c_{\underline{m},\underline{n-1}}+id_{\underline{1}})
\\
&=c_{\underline{m},\underline{n}},
\end{align*}
where in the last equality we make use of the assumption that $\SSS$ is also semistrict.

\bigskip
\noindent
\underline{Proof of (SH3)}. We need to prove that $F_{can}(d'_{[m],[n],[p]})=d'_{\underline{m},\underline{n},\underline{p}}$ for each $m,n,p\geq 0$. When $p=0$ or $p=1$ the result follows from the fact that in every rig category $\SSS$ the right distributors 
$d'_{x,y,0}$ and $d'_{x,y,1}$ are such that the four diagrams
\[
\xymatrix{(x+y)0\ar[r]^{d'_{x,y,0}}\ar[d]_{n_{x+y}} & x0+y0\ar[d]^{n_x+n_y} \\ 0 & 0+0\ar[l]^{l_0}}\quad
\xymatrix{(x+y) 1\ar[rr]^{d'_{x,y,1}}\ar[rd]_{r'_{x+y}}  & & x 1+y 1\ar[ld]^{r'_x+l'_y} \\ & x+y &}
\]
commute for every objects $x,y\in\Ss$. Hence they are identities when $\SSS$ is semistrict and consequently, they are strictly preserved. As to the cases $p\geq 2$, for any given $m,n\geq 0$, they readily follow from (SH2), and Lemma~\ref{lema_distribuidors_dreta} below, which proves that some of the right distributors in every left semistrict rig category are actually determined by the additive commutators.
\qed

\subsubsection{\sc Lemma.}\label{lema_distribuidors_dreta}
Let $\SSS$ be a left semistrict rig category, and let $\underline{n}$ be as before for each $n\geq 0$. Then for every objects $x,y\in\Ss$, and $n\geq 2$ the right distributor $d'_{x,y,\underline{n}}$ is given by
\begin{align}\label{distribuidors_dreta}
d'_{x,y,\underline{n}}&=(id_{x\bullet\underline{n-1}}+c_{y\bullet\underline{n-1},\,x}+id_y)\,(id_{x\bullet\underline{n-2}}+c_{y\bullet\underline{n-2},\,x}+id_{x+y\bullet\underline{2}}) 
\\
&\hspace{2truecm}\cdots (id_{x\bullet\underline{2}}+c_{y\bullet\underline{2},\,x}+id_{x\bullet\underline{n-3}+y\bullet\underline{n-2}})\,(id_{x}+c_{y,\,x}+id_{x\bullet\underline{n-2}+y\bullet\underline{n-1}})\nonumber
\end{align}

\medskip
\noindent
{\sc Proof.} In every rig category $\SSS$ the diagram
\begin{equation} \label{axioma_coherencia_fonamental}
\xymatrix{(x+y)(z+t)\ar[r]^-{d'_{x,y,z+t}}\ar[d]_{d_{x+y,z,t}} & x(z+t)+y(z+t)\ar[r]^{d_{x,z,t}+d_{y,z,t}} & (xz+xt)+(yz+yt)\ar[d]^{v_{xz,xt,yz,yt}}  \\ (x+y)z+(x+y)t\ar[rr]_{\ \ d'_{x,y,z}+d'_{x,y,t}} & & (xz+yz)+(xt+yt)}
\end{equation}
commutes for every objects $x,y,z,t\in\Ss$, where $v_{a,b,c,d}:(a+b)+(c+d)\stackrel{\cong}{\to}(a+c)+(b+d)$ is the canonical isomorphism built from the associator $a$ and the additive commutator $c$ of $\mathscr{S}$. In particular, when $\SSS$ is left semistrict it follows that
\begin{align*}
d'_{x,y,z+t}&=v^{-1}_{xz,xt,yz,yt}\,(d'_{x,y,z}+d'_{x,y,t})
\\
&=(id_{xz}+c_{yz,xt}+id_{yt})\,(d'_{x,y,z}+d'_{x,y,t}).
\end{align*}
Therefore when $z=\underline{n-1}$, $n\geq 2$, and $t=\underline{1}$ we obtain that
\[
d'_{x,y,\underline{n}}=(id_{x\bullet\underline{n-1}}+c_{y\bullet\underline{n-1},x\bullet\underline{1}}+id_{y\bullet\underline{1}})\,(d'_{x,y,\underline{n-1}}+d'_{x,y,\underline{1}}),
\]
a recursive relation whose solution is given by (\ref{distribuidors_dreta}) because, as pointed out before, $d'_{x,y,\underline{1}}=id_{x+y}$ in a semistrict category.
\qed

\subsubsection{\sc Definition.} The strict morphism $\FF_{can}:\widehat{\FF\SSS et}_{sk}\to\SSS$ in Proposition~\ref{morfisme_canonic} is called the {\em canonical homomorphism} of rig categories from $\widehat{\FF\SSS et}_{sk}$ into the left semistrict rig category $\SSS$.

\subsubsection{\sc Remark} What we have proved in Proposition~\ref{morfisme_canonic} is that the commutator of $\SSS$ induces a canonical group homomorphisms $\phi_n:S_n\to Aut_\Ss(\underline{n})$ for each $n\geq 2$, and that the functor from $\widehat{\Ff\Ss et}_{sk}$ to $\Ss$ so defined strictly preserves both the additive and the multiplicative monoidal structures. In fact, the same argument shows the existence of canonical group homomorphisms $\phi_{n;x}:S_n\to Aut_\Ss(nx)$ for anyy object $x$ in $\Ss$ other than the unit object $1$. We just need to replace $c_{1,1}$ by the commutator $c_{x,x}$. Of course, the corresponding functor will no longer extend to a homomorphism of rig categories

\subsubsection{\sc Example} When $\SSS=\widehat{\FF\SSS et}_{sk}$, $\FF_{can}$ is the identity of $\widehat{\FF\SSS et}_{sk}$. In fact, in this case $\underline{n}=[n]$, and $\phi_n$ is the identity.

\subsubsection{\sc Example}
When $\SSS=\MM at_k$ (cf. Example~\ref{Mat_k}), $\FF_{can}$ is the strict morphism with underlying functor the canonical one mapping each object $[n]$ to $n$, and each permutation $\sigma\in S_n$ to the corresponding permutational matrix $P(\sigma)$ obtained from the identity matrix by moving the 1 in the $(i,i)^{th}$-entry to the entry $(\sigma(i),i)$, i.e. the matrix given by $P(\sigma)_{ij}=\delta_{i\,\sigma(j)}$ for each $i,j=1,\ldots,n$. 

\subsubsection{\sc Example}
\label{morfisme_a_endomorfismes_M}
When $\SSS=\EE nd(\mathscr{M})$ for some semistrict symmetric monoidal category $\mathscr{M}=(\Mm,\oplus,{\sf 0})$ (cf. Example~\ref{semianell_endomorfismes}), $\FF_{can}$ is the strict morphism whose underlying functor maps each object $[n]$ to the symmetric monoidal endofunctor
\[
F_{can}[n]\equiv\oplus^n:\mathscr{M}\to\mathscr{M}
\]
given on objects and morphisms by $*\mapsto *\,\oplus\stackrel{n)}{\cdots}\oplus\, *$, and each permutation $\sigma\in S_n$ to the monoidal natural automorphism $F_{can}(\sigma):\oplus^n\Rightarrow\oplus^n$ with $x$-component
\[
F_{can}(\sigma)_x:x\oplus\stackrel{n)}{\cdots}\oplus x\to x\oplus\stackrel{n)}{\cdots}\oplus x
\]
given by the canonical isomorphism defined by the commutator $c_{x,x}$ and the permutation $\sigma$ for each object $x\in\Mm$. This gives the unique (up to equivalence), canonically given $\widehat{\FF\SSS et}$-module structure on every symmetric monoidal category, in the same way as every abelian monoid has a unique, canonicallly given $\NN$-module structure.

\subsection{Main theorem}
We are now able to prove the main result of the paper, namely, that $\widehat{\FF\SSS et}_{sk}$ and hence, $\widehat{\FF\SSS et}$ is biinitial in the 2-category of rig categories. This allows one to view this symmetric 2-rig as the right categorification of the commutative rig $\NN$ of natural numbers.

In order to prove that $\widehat{\FF\SSS et}_{sk}$ is indeed biinitial, we have to prove two things. Firstly, for every (semistrict) rig category $\SSS$ all homomorphisms of rig categories $\FF:\widehat{\FF\SSS et}_{sk}\to\SSS$ must be isomorphic or equivalently, isomorphic to the canonical one $\FF_{can}$. Secondly, for every (semistrict) rig category $\SSS$ the canonical homomorphism $\FF_{can}$ must have the identity as the unique rig transformation to itself. These statements are respectively shown in the next two lemmas.

\subsubsection{\sc Lemma.}\label{lema_principal}
Every homomorphism of rig categories $\FF=(F,\varphi^+,\varepsilon^+,\varphi^{\,\bullet},\varepsilon^{\,\bullet})$ from $\widehat{\FF\SSS et}_{sk}$ to a left semistrict rig category $\SSS$ is isomorphic to the corresponding canonical homomorphism $\FF_{can}$. In fact, the basic maps $(\tau_n)_{n\geq 0}$ of $\FF$ (Definition~\ref{morfismes_basics}) are the components of an invertible rig transformation $\tau:\FF\Rightarrow\FF_{can}$.

\medskip
\noindent
{\sc Proof.} We already know that the basic maps $\tau_n$ are all invertible. Hence we need to prove the following:
\begin{itemize}
\item[(RT1)] $\tau_{n}$ is natural in $[n]$;
\item[(RT2)] $\tau$ is $+$-monoidal;
\item[(RT3)] $\tau$ is $\bullet$\,-monoidal.
\end{itemize}
Since $\FF_{can}$ is strict, (RT2) and (RT3) follow readily from Proposition~\ref{monoidalitat_tau} together with the fact that the first two basic maps are $\tau_0=\varepsilon^+$ and $\tau_1=\varepsilon^\bullet$. In order to prove (RT1), we need to see that the diagram
\begin{equation}\label{naturalitat_xi}
\xymatrix{
F[n]\ar[d]_{F(\sigma)}\ar[r]^{\tau_{n}} & \underline{n}\ar[d]^{F_{can}(\sigma)} 
\\
F[n]\ar[r]_{\tau_{n}} & \underline{n}
}
\end{equation}
commutes for every permutation $\sigma\in S_n$, and every $n\geq 0$. If $n\in\{0,1\}$ this is obvious because $\sigma$ is necessarily an identity. If $n\geq 2$ it is enough to prove it when $\sigma$ is one of the generators $(1,2)_n,\ldots,(n-1,n)_n$ of $S_n$. The generic case follows from the functoriality of $F$ and $F_{can}$ and the invertibility of each $\tau_{n}$.

Since $\tau_{n}$ can be defined recursively by any of the equivalent recurrence relations (\ref{recurrencia_tau_1}) or (\ref{recurrencia_tau_2}), we proceed by induction on $n\geq 2$. If $n=2$ the commutativity of (\ref{naturalitat_xi}) for each $\sigma\in S_2$ just amounts to the commutativity of the outer diagram
\[
\xymatrix{
F[2]\ar[d]_{F((1,2)_{2})}\ar[r]^-{\varphi^+_{[1],[1]}} & F[1]+F[1]\ar[r]^-{\varepsilon^{\,\bullet}+\varepsilon^{\,\bullet}}\ar@{.>}[d]^{c_{F[1],F[1]}} & \underline{2}\ar[d]^{c_{\underline{1},\underline{1}}}\\
F[2]\ar[r]_-{\varphi^+_{[1],[1]}} & F[1]+F[1]\ar[r]_-{\varepsilon^{\,\bullet}+\varepsilon^{\,\bullet}} & \underline{2}
}
\]
where we have already made use of the fact that $F_{can}((1,2)_2)=c_{\underline{1},\underline{1}}$, and the definition of $\tau_{2}$. Now, the right square commutes because of the naturality of $c$, and the left square also commutes because $(1,2)_2=c_{[1],[1]}$ and $F$ is a symmetric $+$-monoidal functor. Let us now assume that for a given $n\geq 2$ the diagram (\ref{naturalitat_xi}) commutes for every generator $\sigma\in S_n$. Then for each $i=1,\ldots,n-1$ the diagram
\[
\xymatrix{
F[n+1]\ar[d]_{F((i,i+1)_{n+1})}\ar[r]^-{\varphi^+_{[n],[1]}} & F[n]+F[1]\ar[r]^-{\tau_{n}+\varepsilon^{\,\bullet}}\ar[d]^{F((i,i+1)_n)+id_{F[1]}} & \underline{n+1}\ar[d]^{F_{can}((i,i+1)_{n+1})} 
\\
F[n+1]\ar[r]_-{\varphi^+_{[n],[1]}} & F[n]+F[1]\ar[r]_-{\tau_{n}+\varepsilon^{\,\bullet}} & \underline{n+1}
}
\]
commutes. Thus the left square is a $\varphi^+$-naturality square, and the right square commutes because
\begin{align*}
F_{can}((i,i+1)_{n+1})&=id_{\underline{i-1}}+c_{\underline{1},\underline{1}}+id_{\underline{n-i}}
\\
&=id_{\underline{i-1}}+c_{\underline{1},\underline{1}}+id_{\underline{n-i-1}}+id_{\underline{1}}
\\
&=F_{can}((i,i+1)_{n})+id_{\underline{1}}
\end{align*}
and by the induction hypothesis. Notice that in the second equality we are using that $i\leq n-1$ and hence, $n-i\geq 1$ to be able to make explicit the identity $id_{\underline{1}}$ on the right. This proves that the diagram
\[
\xymatrix{
F[n+1]\ar[d]_{F(\sigma)}\ar[r]^{\tau_{n+1}} & \underline{n+1}\ar[d]^{F_{can}(\sigma)} 
\\
F[n+1]\ar[r]_{\tau_{n+1}} & \underline{n+1}
}
\]
commutes for all the generators of $S_{n+1}$ except the generator $(n,n+1)_{n+1}$. To prove that the diagram also commutes in this case we make use of the fact that
\[
(n,n+1)_{n+1}=id_{[1]}+(n-1,n)_n
\]
together with the coherence axioms required on $\varphi^+$, from which we know that
\[
(id_{F[1]}+\varphi^+_{[n-1],[1]})\,\varphi^+_{[1],[n]}=(\varphi^+_{[n-1],[1]}+id_{F[1]})\,\varphi^+_{[n],[1]}.
\]
It follows that the previous square with $\sigma=(n,n+1)_{n+1}$ looks like
\[
\xymatrix{
F[n+1]\ar[d]_{F(id_{[1]}+(n-1,n)_{n})}\ar[r]^-{\varphi^+_{[1],[n]}} & F[1]+F[n]\ar[d]^{id_{F[1]}+F((n-1,n)_n)}\ar[r]^{id_{F[1]}+\varphi^+_{[n-1],[1]}\hspace{1truecm}} & F[1]+F[n-1]+F[1]\ar[r]^-{(\varphi_{[1],[n-1]}+id_{F[1]})^{-1}} & F[n]+F[1]\ar[r]^-{\tau_{n}+\varepsilon^{\,\bullet}} & \underline{n+1}\ar[d]_{F_{can}((n,n+1)_{n+1})} 
\\
F[n+1]\ar[r]_-{\varphi^+_{[1],[n]}} & F[1]+F[n]\ar[r]_-{id_{F[1]}+\varphi^+_{[n-1],[1]}\hspace{1truecm}} & F[1]+F[n-1]+F[1]\ar[r]_-{(\varphi_{[1],[n-1]}+id_{F[1]})^{-1}} & F[n]+F[1]\ar[r]_-{\tau_{n}+\varepsilon^{\,\bullet}} & \underline{n+1}
}
\]
The left square commutes because of the naturality of $\varphi^+$. As to the right part, notice that
\[
F_{can}((n,n+1)_{n+1})=id_{\underline{n-1}}+c_{\underline{1},\underline{1}}=id_{\underline{1}}+id_{\underline{n-2}}+c_{\underline{1},\underline{1}}=id_{\underline{1}}+F_{can}((n-1,n)_n)
\]
because $n\geq 2$. Moreover
\[
(\tau_{n}+\varepsilon^{\,\bullet})\,(\varphi_{[1],[n-1]}+id_{F[1]})^{-1}=\tau_{n}\,(\varphi_{[1],[n-1]})^{-1}+\varepsilon^{\,\bullet}=\varepsilon^{\,\bullet}+\tau_{n-1}+\varepsilon^{\,\bullet}.
\]
Therefore we have
\begin{align*}
(\tau_{n}+\varepsilon^{\,\bullet})\,(\varphi_{[1],[n-1]}+id_{F[1]})^{-1}\,(id_{F[1]}+\varphi^+_{[n-1],[1]})&=(\varepsilon^{\,\bullet}+\tau_{n-1}+\varepsilon^{\,\bullet})\,(id_{F[1]}+\varphi^+_{[n-1],[1]})
\\
&=\varepsilon^{\,\bullet}+\tau_{n},
\end{align*}
so that the right subdiagram in the above diagram is nothing but the diagram
\[
\xymatrix{
F[n]+F[1]\ar[d]_{id_{F[1]}+F((n-1,n)_n)}\ar[r]^-{\varepsilon^{\,\bullet}+\tau_{n}} & \underline{n+1}\ar[d]^{id_{\underline{1}}+F_{can}((n-1,n)_n)} 
\\
F[n]+F[1]\ar[r]_-{\varepsilon^{\,\bullet}+\tau_{n}} & \underline{n+1}
}
\]
whose commutativity follows from the functoriality of $+$ and the induction hypothesis.
\qed

\subsubsection{\sc Lemma}\label{unicitat} $\FF_{can}:\widehat{\FF\SSS et}_{sk}\to\SSS$ has the identity as unique (2-)endomorphism.

\medskip
\noindent
{\sc Proof.} Since $\SSS$ is assumed to be left semistrict and $\FF_{can}$ is strict, it follows from Definition~\ref{transformacio_rig} that an endomorphism $\xi:\FF_{can}\Rightarrow\FF_{can}$ consists of a collection of morphisms $\xi_n:\underline{n}\to\underline{n}$ in $\Ss$, for $n\geq 0$, such that the following diagrams commute: 
\begin{itemize}
\item[(1)] (naturality) for each $n\geq 0$ and each permutation $\sigma\in S_n$
\[
\xymatrix{
\underline{n}\ar[r]^{\xi_{n}}\ar[d]_{F_{can}(\sigma)} & \underline{n}\ar[d]^{F_{can}(\sigma)}
\\
\underline{n}\ar[r]_{\xi_{n}} & \underline{n}\,;
}
\]
\item[(2)] ($+$-monoidality) $\xi_0=id_{\underline{0}}$, and for each $m,n\geq 0$
\[
\xymatrix{
\underline{m+n}\ar[r]^{\xi_{m+n}}\ar[d]_{id} & \underline{m+n}\ar[d]^{id}
\\
\underline{m+n}\ar[r]_{\xi_{m}+\xi_{n}} & \underline{m+n}\,;
}
\] 
\item[(3)] ($\bullet$-monoidality) $\xi_1=id_{\underline{1}}$, and for each $m,n\geq 0$
\[
\xymatrix{
\underline{mn}\ar[r]^{\xi_{mn}}\ar[d]_{id} & \underline{mn}\ar[d]^{id}
\\
\underline{mn}\ar[r]_{\xi_{m}\bullet\xi_{n}} & \underline{mn}\,.
}
\] 
\end{itemize}
In particular, it must be $\xi_{n+1}=\xi_n+\xi_1$ for each $n\geq 0$. Together with the fact that $\xi_0,\xi_1$ are identities, it follows that $\xi_n=id_{\underline{n}}$ for each $n\geq 0$.
\qed

\medskip
Together, these two lemmas imply that there exists one and only one isomorphism between every two homomorphisms $\FF,\FF':\widehat{\FF\SSS et}_{sk}\to\SSS$ as above, necessarily given by the composite
\[
\xymatrix{
\FF\ar@{=>}[r]^\tau & \FF_{can}\ar@{=>}[r]^{(\tau')^{-1}} & \FF'\,}
\]
with $\tau,\tau'$ the invertible rig transformations in Lemma~\ref{lema_principal} defined by the respective basic maps of $\FF$ and $\FF'$. Thus we have proved the following.

\subsubsection{\sc Theorem}\label{main} The symmetric 2-rig $\widehat{\FF\SSS et}$ of finite sets is biinitial in the 2-category $\mathbf{RigCat}$.

\subsection{Explicit description of the homomorphisms from the semistrict 2-rig of finite sets}\label{conclusions}

Let us denote by $\Hh om_{\mathbf{RigCat}}(\widehat{\FF\SSS et}_{sk},\SSS)$ the category of rig category homomorphisms from $\widehat{\FF\SSS et}_{sk}$ to $\SSS$. According to Theorem~\ref{main}, it is a contractible category for any $\SSS$. In fact, when $\SSS$ is left semistrict we have proved more than that, in the sense that we have a quite explicit description of the objects in this category (and of the unique isomorphism between any of them).

Thus it follows from the commutativity of (\ref{naturalitat_xi}) that the action on morphisms of the underlying functor $F$ of each $\FF:\widehat{\FF\SSS et}_{sk}\to\SSS$ is determined by the basic maps $\tau_n:F[n]\to\underline{n}$, $n\geq 0$, of $\FF$. Since the data $(\varphi^+,\varepsilon^+,\varphi^\bullet,\varepsilon^\bullet)$ is also determined by these maps (cf. Proposition~\ref{monoidalitat_tau}), we conclude that every homomorphism $\FF$ is completely given by its action on objects, and the basic maps. In other words, $\FF$ is given by just two sequences, a sequence of objects $x=(x_n)_{n\geq 0}$ in $\Ss$ such that $x_n\cong\underline{n}$ for each $n\geq 0$, and a sequence of isomorphisms $\tau=(\tau_n)_{n\geq 0}$, with $\tau_n:x_n\to\underline{n}$. It turns out that these sequences can be chosen arbitrarily. More precisely, we have the following.

\subsubsection{\sc Proposition}\label{morfismes_normalitzats}
For any sequence $x=(x_n)_{n\geq 0}$ of objects in $\Ss$, with $x_n\cong\underline{n}$, and any sequence of isomorphisms $\tau=(\tau_n)_{n\geq 0}$, with $\tau_n:x_n\to\underline{n}$, the functor $F(x,\tau):\widehat{\Ff\Ss et}_{sk}\to\Ss$ defined on objects and morphisms respectively by
\begin{align}
    F(x,\tau)[n]&=x_n,\quad n\geq 0,\label{F_tau_0}
    \\
    F(x,\tau)(\sigma)&=\tau^{-1}_n\,F_{can}(\sigma)\,\tau_n,\quad \sigma\in S_n,\label{F_tau_1}
\end{align}
together with the isomorphisms
\begin{align}
    \varphi(\tau)^+_{[m],[n]}&=(\tau_m+\tau_n)^{-1}\,\tau_{m+n},\quad m,n\geq 0, \label{varphi_tau_+}
    \\
    \varphi(\tau)^\bullet_{[m],[n]}&=(\tau_m\bullet\tau_n)^{-1}\tau_{mn},\quad m,n\geq 0, \label{varphi_tau_bullet}
\end{align}
and $\varepsilon(\tau)^+=\tau_0$, $\varepsilon(\tau)^\bullet=\tau_1$ is a homomorphism of rig categories $\FF(x,\tau)$. 

\medskip
\noindent
{\sc Proof.} Since both rig categories are left semistrict, the required naturality and coherence conditions on the data $(\varphi(\tau)^+,\varepsilon(\tau)^+,\varphi(\tau)^\bullet,\varepsilon(\tau)^\bullet)$ reduce to the following:
\begin{itemize}
\item[(A1)] ({\em naturality of} $\varphi(\tau)^+$) the diagram
\[
\xymatrix{
x_{m+n}\ar[r]^{\varphi(\tau)^+_{[m],[n]}}\ar[d]_{F(x,\tau)(\rho+\sigma)} & x_m+x_n\ar[d]^{F(x,\tau)(\rho)+F(x,\tau)(\sigma)}
\\
x_{m+n}\ar[r]_{\varphi(\tau)^+_{[m],[n]}} & x_m+x_n
}
\]
commutes for each permutations $\rho\in S_m$ and $\sigma\in S_n$, and each $m,n\geq 0$;

\item[(A2)] ({\em coherence axiom on} $\varphi(\tau)^+$) the diagram
\[
\xymatrix{
x_{m+n+p}\ar[rr]^{\varphi(\tau)^+_{[m],[n+p]}}\ar[d]_{\varphi(\tau)^+_{[m+n],[p]}} && x_m+x_{n+p}\ar[d]^{id_{x_m}+\varphi(\tau)^+_{[n],[p]}}
\\
x_{m+n}+x_p\ar[rr]_{\varphi(\tau)^+_{[m],[n]}+id_{x_p}} && x_m+x_n+x_p
}
\]
commutes for each $m,n,p\geq 0$;

\item[(A3)] ({\em coherence axiom on} $\varphi(\tau)^+$) the diagram
\[
\xymatrix{
x_{m+n}\ar[r]^{\varphi(\tau)^+_{[m],[n]}}\ar[d]_{F(x,\tau)(c_{[m],[n]})} & x_m+x_n\ar[d]^{c_{x_m,x_n}}
\\
x_{n+m}\ar[r]_{\varphi(\tau)^+_{[n],[m]}} & x_n+x_m
}
\]
commutes for each $m,n\geq 0$;

\item[(A4)] ({\em coherence axiom on} $\varepsilon(\tau)^+$) the diagrams
\[
\xymatrix{
x_n\ar[r]^-{\varphi(\tau)^+_{[0],[n]}}\ar[dr]_{id_{x_n}} & x_0+x_n\ar[d]^{\varepsilon(\tau)^++id_{x_n}} 
\\
& x_n}
\quad 
\xymatrix{
x_n\ar[r]^-{\varphi(\tau)^+_{[n],[0]}}\ar[dr]_{id_{x_n}} & x_n+x_0\ar[d]^{id_{x_n}+\varepsilon(\tau)^+} 
\\
& x_n}
\]
commute for each $n\geq 0$;

\item[(A5)] ({\em naturality of} $\varphi(\tau)^\bullet$) the diagram
\[
\xymatrix{
x_{mn}\ar[r]^-{\varphi(\tau)^\bullet_{[m],[n]}}\ar[d]_{F(x,\tau)(\rho\bullet\sigma)} & x_m\bullet x_n\ar[d]^{F(x,\tau)(\rho)\bullet F(x,\tau)(\sigma)}
\\
x_{mn}\ar[r]_{\varphi(\tau)^\bullet_{[m],[n]}} & x_m\bullet x_n
}
\]
commutes for each permutations $\rho\in S_m$ and $\sigma\in S_n$, and each $m,n\geq 0$;

\item[(A6)] ({\em coherence axiom on} $\varphi(\tau)^\bullet$) the diagram
\[
\xymatrix{
x_{mnp}\ar[rr]^{\varphi(\tau)^\bullet_{[m],[np]}}\ar[d]_{\varphi(\tau)^\bullet_{[mn],[p]}} && x_m\bullet x_{np}\ar[d]^{id_{x_m}\bullet\varphi(\tau)^\bullet_{[n],[p]}}
\\
x_{mn}\bullet x_p\ar[rr]_{\varphi(\tau)^\bullet_{[m],[n]}\bullet id_{x_p}} && x_m\bullet x_n\bullet x_p
}
\]
commutes for each $m,n,p\geq 0$;

\item[(A7)] ({\em coherence axiom on} $\varepsilon(\tau)^\bullet$) the diagrams
\[
\xymatrix{
x_n\ar[r]^-{\varphi(\tau)^\bullet_{[1],[n]}}\ar[dr]_{id_{x_n}} & x_1\bullet x_n\ar[d]^{\varepsilon(\tau)^\bullet\bullet id_{x_n}} 
\\
& x_n}
\quad 
\xymatrix{
x_n\ar[r]^-{\varphi(\tau)^\bullet_{[n],[1]}}\ar[dr]_{id_{x_n}} & x_n\bullet x_1\ar[d]^{id_{x_n}\bullet\varepsilon(\tau)^\bullet} 
\\
& x_n}
\]
commute for each $n\geq 0$;

\item[(A8)] ({\em coherence axiom between $\varphi(\tau)^+$ and} $\varphi(\tau)^\bullet$) the diagram
\[
\xymatrix{
x_{mn+mp}\ar[d]_{\varphi(\tau)^{\,\bullet}_{[m],[n+p]}}\ar[rr]^-{\varphi(\tau)^+_{[mn],[mp]}} && x_{mn}+x_{mp}\ar[d]^{\varphi(\tau)^{\,\bullet}_{[m],[n]}\,+\,\varphi(\tau)^{\,\bullet}_{[m],[p]}} 
\\
x_m\bullet x_{n+p}\ar[rr]_-{id_{x_m}\bullet\, \varphi(\tau)^+_{[n],[p]}} && x_m\bullet x_n+x_m\bullet x_p;
}
\]
commutes for each $m,n,p\geq 0$;

\item[(A9)] ({\em coherence axiom between $\varphi(\tau)^+$ and} $\varphi(\tau)^\bullet$) the diagram
\[
\xymatrix{
x_{(m+n)p}\ar[d]_{\varphi(\tau)^{\,\bullet}_{[m+n],[p]}}\ar[rr]^-{F(x,\tau)(d'_{[m],[n],[p]})} && x_{mp+np}\ar[rr]^-{\varphi(\tau)^+_{[mp],[np]}} && x_{mp}+x_{np}\ar[d]^{\varphi(\tau)^{\bullet}_{[m],[p]}\,+\,\varphi(\tau)^{\bullet}_{[n],[p]}}
\\
x_{m+n}\bullet x_p\ar[rr]_-{\varphi(\tau)^+_{[m],[n]}\bullet\, id_{x_p}} && (x_m+x_n)\bullet x_p\ar[rr]_-{d'_{x_m,x_n,x_p}} && x_m\bullet x_p+x_n\bullet x_p
} 
\]
commutes for each $m,n,p\geq 0$.
\end{itemize}
When made explicit using (\ref{F_tau_1})-(\ref{varphi_tau_bullet}), all conditions are easy to check. For instance, condition (A1) (and condition (A5) is similar but with $\bullet$ instead of $+$) amounts to 
\[
[\tau_m^{-1}\,F_{can}(\rho)\,\tau_m+\tau_n^{-1}\,F_{can}(\sigma)\,\tau_n]\,(\tau_m+\tau_n)^{-1}\,\tau_{m+n}=(\tau_m+\tau_n)^{-1}\,\tau_{m+n}\,\tau_{m+n}^{-1}\,F_{can}(\rho+\sigma)\,\tau_{m+n},
\]
and this equality is a consequence of the functoriality of $+$ (resp. $\bullet$) together with the fact that $F_{can}$ preserves the sum (resp. the product) of morphisms. (A2), and its analog (A6), simply follow from the functoriality of $+$ and $\bullet$, respectively. (A3) is a consequence of the functoriality of $+$, the naturality of $c_{x_m,x_n}$ in $x_m,x_n$, and the fact that $F_{can}$ preserves the commutators. (A4), and its analog (A7), follow from the fact that the left and right additive and multiplicative unitors of $\SSS$ are assumed to be identities. As to (A8), it is a consequence of the functoriality of $+$ and $\bullet$ together with the fact that the left distributor of $\SSS$ is assumed to be trivial. Finally, (A9) follows from the functoriality of $+$ and $\bullet$, the naturality of $d'_{x_m,x_n,x_p}$ in $x_m,x_n,x_p$, and the fact that $F_{can}$ preserves the right distributors. 
\qed

\subsubsection{\sc Corollary}
For every left semistrict rig category $\SSS$ the set of objects in $\Hh om_{\mathbf{RigCat}}(\widehat{\FF\SSS et}_{sk},\SSS)$ is in bijection with the set of pairs $(x,\tau)$, with $x=(x_n)_{n\geq 0}$ any sequence of objects in $\Ss$ such that $x_n\cong\underline{n}$ for each $n\geq 0$, and $\tau=(\tau_n)_{n\geq 0}$ any sequence of isomorphisms $\tau_n:x_n\to\underline{n}$. In particular, if $\SSS$ is skeletal $\Hh om_{\mathbf{RigCat}}(\widehat{\FF\SSS et}_{sk},\SSS)$ is isomorphic to the contractible category with set of objects $\prod_{n\geq 0}\,Aut_\Ss(\underline{n})$.

\subsubsection{\sc Example} Let $\SSS$ be the symmetric 2-rig $\widehat{\FF\SSS et}_{sk}$ itself. It follows that $\Ee nd_{\mathbf{RigCat}}(\widehat{\FF\SSS et}_{sk})$ is isomorphic to the contractible category with $\prod_{n\geq 0}\,S_n$ as set of objects. In fact, as the category of endomorphisms of any object in a 2-category, it is a strict monoidal category with the tensor product given by the composition of endomorphisms. Even more, since every endomorphism of $\widehat{\FF\SSS et}_{sk}$ is strictly invertible, $\Ee nd_{\mathbf{RigCat}}(\widehat{\FF\SSS et}_{sk})$ is a {\em strict 2-group}, i.e. a strict monoidal groupoid all of whose objects are strictly invertible. Then it is easy to check that the above bijection between its set of objects, equipped with the group law induced by the tensor product, and $\prod_{n\geq 0}\,S_n$ with its usual direct product group structure is a group isomorphism.

\bibliography{conjectura_baez}{}
\bibliographystyle{plain}

\end{document}